\documentclass[11pt,a4paper]{article}
\usepackage{graphicx} 
\usepackage{pifont}
\usepackage{bigints}
\usepackage{float}
\usepackage[margin=1.70cm]{geometry}
\usepackage{authblk}
\usepackage{makecell}   
 
 \usepackage{setspace}
 \usepackage{lscape}
 \usepackage{array}
 \usepackage{tikz}
\usetikzlibrary{decorations.pathreplacing}
\usepackage{booktabs}
 \usepackage{multirow}
 \usepackage{amsmath}
 \usepackage{amsthm}
 \usepackage[utf8]{inputenc}
\usepackage[pdfencoding=auto]{hyperref}
 \usepackage{longtable}
 \usepackage{amsfonts}
\usepackage{booktabs}
\usepackage{siunitx}
\newcolumntype{P}[1]{>{\raggedright\arraybackslash}p{#1}}
 \newtheorem{theorem}{Theorem}[section] 
 \newtheorem{proposition}{Proposition} [section]
 \newtheorem{definition}{Definition}[section]
 \newtheorem{example}{Example}[section]
 \newtheorem{remark}{Remark}[section]

 \usepackage{tabularx}
\usepackage{amssymb}
 \usepackage{caption}
 \usepackage{subcaption}
 \usepackage{array}
 \usepackage{xcolor}
 \usepackage{hyperref}
 \hypersetup{
    colorlinks=true,
    linkcolor=blue,
    filecolor=blue,      
    urlcolor=blue,
    citecolor=blue
}
 \usepackage[T1]{fontenc}
 \newcommand{\R}{\ensuremath{\mathbb{R}}}

\begin{document}

     \title{\bf On Generalized Likelihood Estimation Based on the Logarithmic Norm Relative Entropy}
     \date{}
 \author[1]{Himanshi Singh}
 \author[2]{Abhik Ghosh}
 \author[1]{Nil Kamal Hazra}
\affil[1]{Department of Mathematics, Indian Institute of Technology Jodhpur, Rajasthan, India}
\affil[2]{Interdisciplinary Statistical Research Unit, Indian
Statistical Institute, Kolkata, India}
   \maketitle
\begin{abstract}
Traditional likelihood based methods for parameter estimation get highly affected when the given data is contaminated by outliers even in a small proportion. 
In this paper, we consider a robust parameter estimation method, namely the minimum logarithmic norm relative entropy (LNRE) estimation procedure, and study different (generalized) sufficiency principles associated with it. We introduce a new two-parameter power-law family of distributions (namely, $\mathcal{M}^{(\alpha,\beta)}$-family), which is shown to have a fixed number of sufficient statistics, independent of the sample size, with respect to the generalized likelihood function associated with the LNRE. Then, we obtain the generalized minimal sufficient statistic for this family and derive the generalized Rao-Blackwell theorem and the generalized Cram\'{e}r-Rao lower bound for the minimum LNRE estimation. We also study the minimum LNRE estimators (MLNREEs) for the family of Student's distributions particularly in detail. Our general results reduces to the classical likelihood based results under the exponential family of distributions at specific choices of the tuning parameter $\alpha$ and $\beta$.  Finally, we present simulation studies followed by a real data analysis, which highlight the practical utility of the MLNREEs for data contaminated by possible outliers. Along the way we also correct a mistake found in a recent paper  on related theory of generalized likelihoods. 
\end{abstract}



\section{Introduction}\label{S1}
Divergence (or distance) is a real-valued function used to quantify the dissimilarity between a pair of probability distributions. 
There are many popular (density-based) divergence measures, namely, the Kullback-Leibler divergence (KLD), the Renyi divergence, the Tsallis divergence (also commonly known as the power divergence (PD)), the density power divergence (DPD), the relative $\alpha$-entropy (or the logarithmic density power divergence (LDPD)), the logarithmic norm relative entropy (LNRE), etc. (see \cite{DPD,Basu_book,LNRE,LDPD,KLD,Entropy_Review,Renyi,Tsallis,Divergence_review}). Minimum divergence  methods are widely used in robust statistical inference.
Minimization problems associated with the KLD gave rise the classical notion of information projection  (\cite{Csiszar,KLD}). In the literature, there are two types of information projections, namely, forward projection (i.e., minimization of a divergence with respect to the first argument) and reverse projection (i.e., minimization of a divergence  with respect to the second argument). 
In case of forward projection, the KLD behaves like the ``squared Euclidean distance'' and satisfies the ``Pythagorean property'', whereas  
reverse projection of the KLD is related with the classical maximum likelihood (ML) estimation. 
\\\hspace*{0.2 in}
The estimation of parameter(s) for a given parametric family of probability distributions (often called models) is highly important which is solved by suitable estimators, defined as appropriate functions of random sample. In particular, the estimators based on sufficient statistics always play a key role in statistical decision theory (\cite{Sufficient}). 
Note that the entire sample is always a sufficient statistic but it is not helpful in data reduction. A sufficient statistic is useful when its dimension is fixed, independent of the sample size. This problem was addressed in Fisher-Darmois-Koopman-Pitman theorem (\cite{Darmois,Fisher, Koopman, Pitman}), which states that, under certain regularity conditions, the exponential family is the only one that has a fixed number of sufficient statistics, independent of the sample size. One may note that this statement may not be true if any other estimation procedure apart from the ML estimation is under consideration.  
The ML estimation (i.e., maximization of the usual likelihood function) has a close connection with the minimum distance estimation based on the KLD under the parametric estimation problem  as follows. Let $ Y_1^n := \{Y_1, \dots , Y_n\}$ be an independent and identically distributed (i.i.d.) random sample drawn according to a distribution having the probability density function (pdf) $g$, and let $y_1^n := \{y_1,\dots,y_n\}$ be the realized value of $Y_1^n$. Let $f_\lambda \in \mathcal{P}: = \{f_\lambda : \lambda \in \Lambda \subset \R^k\} $ be the model pdf with parameter $\lambda \in \Lambda$, the parameter space. Then, obtaining the maximum likelihood estimate (MLE) of the parameter is equivalent to solve the minimization of the KLD for the pair $(\widehat{g}, f_\lambda)$ with respect to $\lambda$, where $\widehat{g}$ is a suitable estimate of $g$ based on an observed realization $y_1^n$. 
However, the MLE is strongly affected by outliers and hence, it is not appropriate to use for contaminated data. Thus, several distance-based robust estimation methods  were introduced in the literature (see \cite{DPD, LNRE, LDPD, LSD_a, LSD}). 
Among the existing such methods,  the minimum DPD estimation and the minimum LDPD estimation are two popular ones. Similar to the ML estimation, these estimations are also associated with some generalized likelihood functions.  For example, the minimum DPD estimator (MDPDE) can be obtained by maximizing the Basu et al. likelihood function, whereas the minimum LDPD estimator (MLDPDE) is obtained by maximizing the  Jones et al. likelihood function.
Based on these generalized likelihood functions, Gayen and Kumar \cite{Generalized_sufficiency} recently introduced the notion of generalized sufficiency and studied different sufficiency principles associated with the MDPDE and MLDPDE. 
%
\\\hspace*{0.2 in}Like DPD and LDPD, the Logarithmic norm relative entropy (LNRE) is another useful family of divergence measures which was introduced by Ghosh and Basu~\cite{LNRE} as a two-parameter generalization of the LDPD. The key importance of this divergence and its usefulness in robust parameter estimation were well described in  \cite{LNRE, Scale_invariant, LSD_a,LSD}.  The LNRE is {a reformulation} of the logarithmic super divergence (LSD) {with an extended range of tuning parameters} (\cite{LSD}) and contains many important divergences (e.g., KL, LDPD, etc.) as special cases. Statistical inferences based on the LNRE provide a better trade-off between robustness and efficiency compared to those based on other related divergences 
(\cite{LSD_a, LSD}). The minimum LNRE estimator (MLNREE) is obtained by minimizing LNRE or equivalently, maximizing the generalized likelihood function associated with the LNRE (derived below in Section \ref{S7}). Different asymptotic properties of the MLNREE were studied by  Maji et al.~\cite{LSD}. Further, the MLNREE for different discrete distributions (namely, Poisson, Geometric and Binomial) were obtained in \cite{LNRE, Scale_invariant,LSD}. They also studied the robustness of these estimators by performing several simulation studies. However, a similar study for any continuous distributions is yet to be explored. 
\\\hspace*{0.2 in} Keeping the importance of the MLNREE in robust statistical inference in mind, we focus to study different sufficiency principles associated with the minimum LNRE estimation procedure. Moreover, we derive the MLNREE for the family of Student's distributions (see Definition \ref{s_d}) and show that this performs better than  MDPDE and MLDPDE for sample data contaminated by outliers.
To be more specific, what we consider here and what mostly constitutes the novelty and contribution of our paper, is as follows:

\begin{itemize}
    \item We introduce a new two-parameter power-law family of distributions, namely $\mathcal{M}^{(\alpha, \beta)}$-family. Then, we show that this family has a fixed number of sufficient statistics, independent of sample size, with respect to the generalized likelihood function associated  with the LNRE. Further, we obtain the generalized minimal sufficient statistic for regular $\mathcal{M}^{(\alpha, \beta)}$-family.
    \item We obtain the generalized Rao-Blackwell theorem for the minimum LNRE estimation. 
   Further, we derive the generalized Cram\'{e}r-Rao lower bound for MLNREEs of $\mathcal{M}^{(\alpha, \beta)}$- family. We show that the generalized minimal sufficient statistic of regular $\mathcal{M}^{(\alpha, \beta)}$- family does not attain this lower bound, apart from a particular case.
    \item In the process of developing the above-mentioned results, we correct a mistake found in a recent paper. In Theorem $20$ of Gayen and Kumar \cite{Generalized_sufficiency}, it was claimed that the generalized minimal sufficient statistic of regular $\mathcal{M^{(\alpha)}}$-family is the best estimator of its expected value. However, we show that this claim is not true in general by providing a counterexample.
     \item We obtain the MLNREEs for the family of Student's distributions with both positive and negative degrees of freedom. Further, by suitable simulation studies and a real data analysis, we demonstrate that the MLNREEs performs better than the MDPDE and the MLDPDE for data contaminated by outliers.
\end{itemize}
\hspace*{0.2 in}The rest of the paper is organized as follows. In Section \ref{S2}, we present the required preliminaries for the sake of completeness. In Section \ref{S3}, we derive the estimating equations for the minimum LNRE estimation and define the generalized likelihood function associated with the LNRE. In Section \ref{S4}, we introduce the $\mathcal{M}^{(\alpha, \beta)}$-family and derive the generalized Fisher-Darmois-Koopman-Pitman theorem for this family.  Further, we obtain the minimal sufficient statistic for this family. 
In Section \ref{S5}, we discuss the generalized Rao-Blackwell theorem for the minimum LNRE estimation. In Section \ref{S6}, we derive the generalized Cram\'{e}r-Rao lower bound for the minimum LNRE estimation and verify whether the generalized minimal sufficient statistic of regular $\mathcal{M}^{(\alpha, \beta)}$-family attains this bound. In Section \ref{S7}, we derive the MLNREEs for the family of Student's distributions and present their empirical performances. 
Finally, some concluding remarks are given in Section \ref{S8}. A detailed analysis illustrating the computation of the MLNREEs for the family of Student's $r$ distributions (having  negative degrees of freedom) is provided in Appendix.
\section{Preliminaries}\label{S2}
We begin this section by introducing required notation followed by necessary definitions and preliminary results. We denote the set of real numbers and the set of positive real numbers by $\mathbb{R}$ and $\mathbb{R}^{+}$, respectively. 
 In what follows, we give the formal definitions of some divergence measures mentioned in the introduction (Section \ref{S1}) that will be used throughout this paper.
\begin{definition}
Let $f$ and $g$ be two probability density functions having the common support $S\subset \mathbb{R}$. Further, let $\alpha > 0$ and $\beta \in \R$ be two tuning parameters.
\begin{enumerate}
    \item [$(i)$] The Kullback-Leibler divergence (KLD) (also known as the relative entropy), denoted  by $ \mathcal{D}_{KL} (g,f)$, is given by \cite{KLD}
\begin{equation}\label{E1}
  \mathcal{D}_{KL} (g,f)=  \int g(y) \log \frac{g(y)}{f(y)} dy. 
\end{equation}
\item [$(ii)$] The density power divergence (DPD), denoted by $\mathcal{B}_\alpha (g,f)$, is given by \cite{DPD} 
\begin{equation}\label{E2}
  \mathcal{B}_\alpha (g,f)= \frac{\alpha}{1-\alpha} \int g(y)f^{\alpha-1}(y)dy - \frac{1}{1-\alpha} \int g^\alpha(y) dy + \int f^\alpha(y) dy, \quad \alpha \neq 1.
\end{equation}
The case $\alpha=1$ is defined in a limiting sense (as $\alpha\rightarrow 1$) and coincides with the KLD.
\item [$(iii)$] The logarithmic density power divergence (LDPD), denoted by $\mathcal{I}_\alpha (g,f)$, is given by \cite{LDPD}
\begin{equation}\label{E3}
\mathcal{I}_\alpha (g,f)= \frac{\alpha}{1-\alpha} \log \int g(y)f^{\alpha-1}(y)dy - \frac{1}{1-\alpha} \log \int g^\alpha (y)dy + \log \int f^\alpha(y) dy, \quad \alpha \neq 1.
\end{equation}
The case $\alpha=1$ is again defined in a limiting sense (as $\alpha\rightarrow 1$) and coincides with the KLD.
\item [$(iv)$] The logarithmic norm relative entropy (LNRE), denoted by $\mathcal{RE}_{\alpha,\beta}^{\mathcal{LN}} (g,f)$, is given by \cite{LNRE}
\begin{equation}\label{E4}
\mathcal{RE}_{\alpha,\beta}^{\mathcal{LN}} (g,f) = \frac{\alpha}{ \beta(\beta-\alpha)} \log \int g^\beta(y) f^{\alpha-\beta}(y)dy - \frac{1}{\beta-\alpha} \log \int g^\alpha(y) dy + \frac{1}{\beta}\log \int f^\alpha(y) dy, \quad \alpha \neq \beta.     
\end{equation}
The case $\alpha=\beta$ is defined in the limiting sense as $\beta \to \alpha$  and resulting divergences has the form $$\mathcal{RE^{LN}_{\alpha,\alpha}}(g,f)=\frac{\int g^\alpha\log\left(\frac{g}{f}\right)dy}{\int g^\alpha dy} + \frac{1}{\alpha} \log\left(\frac{\int f^\alpha dy}{\int g^\alpha dy}\right).$$
Moreover, \eqref{E4} coincides with the KLD when $\alpha\rightarrow1$ and $\beta\rightarrow1$. $\hfill\Box$
\end{enumerate}
\end{definition}
Let us now consider the problem of parametric estimation as mentioned in Section \ref{S1}, with $\lambda$ being the $k$-dimensional parameter of interest. It can be estimated by minimizing $D(\widehat{g}, f_\lambda)$ with respect to $\lambda$, where $\widehat{g}$ is an empirical estimate of $g$ based on the sample and $D$ is any suitably defined divergence measure. In particular, the resulting estimating equations corresponding to the divergences given in \eqref{E1}, \eqref{E2} and \eqref{E3} turns out to be
\begin{equation}\label{E5}
   \frac{1}{n} \sum_{i=1}^n u(y_j, \lambda) = 0,
\end{equation}

\begin{equation}\label{E6}
\frac{1}{n} \sum_{j=1} ^ {n} f_\lambda^{\alpha-1}(y_j)u(y_j,\lambda) = \int f_\lambda^\alpha(y)u(y;\lambda)  dy,
\end{equation}
and
\begin{equation}\label{E7}
\frac{\frac{1}{n} \sum_{j=1} ^ {n}f_\lambda^{\alpha-1}(y_j)u(y_j;\lambda)}{\frac{1}{n} \sum_{j=1} ^ {n} f_\lambda^{\alpha-1}(y_j)} = \frac{\int f_\lambda^\alpha(y)u(y;\lambda)  dy}{\int f_\lambda^\alpha(y)  dy},
\end{equation}    
respectively, where $u(y; \lambda)= \nabla \log ~ f_\lambda(y)$, and $\nabla$ denotes the gradient with respect to $\lambda$ (see \cite{DPD, LDPD, Windham}). 
\\\hspace*{0.2 in}
Recall that maximization of a generalized likelihood function associated with a divergence with respect to parameter of interest is the same as that of minimization of the underlying divergence. Accordingly, generalized likelihood functions associated with the divergences given in \eqref{E1}, \eqref{E2} and \eqref{E3} are given by \cite{Generalized_sufficiency}

\begin{equation}\label{E8}
 \mathcal{L}(y_1^n, \lambda) :=  \frac{1}{n}\sum_{j=1}^n \log \hspace{0.05 in} f_\lambda (y_j),     
\end{equation}

\begin{equation}\label{E9}
 \mathcal{L}_{\mathcal{B}}^{(\alpha)}(y_1 ^n;\lambda):=\frac{1}{n}\sum_{j=1} ^ {n} \left[\frac{ \alpha f_\lambda^{\alpha-1}(y_j) - 1}{\alpha-1}\right] -  \int f_\lambda^\alpha(y)  dy, 
\end{equation}
and
\begin{equation}\label{E10}
 \mathcal{L}_{\mathcal{J}}^{(\alpha)}(y_1 ^n;\lambda) :=\frac{1}{\alpha-1}\log \left[\frac{1}{n}\sum_{j=1} ^ {n} f_\lambda^{\alpha-1}(y_j)\right] - \frac{1}{\alpha}\log \left[\int f_\lambda^\alpha(y)  dy\right],
\end{equation}
respectively. In what follows, we give the definition of generalized sufficiency and discuss some associated notions following \cite{Generalized_sufficiency}. In this respect, let us again consider the set-up and notation of the parametric estimation problem as mentioned in Section \ref{S1}. Let $T := T(Y_1^n)$ be a statistic (which can be a vector-valued function) and let $\mathcal{T} = \{ T(x_1^n): x_1^n \in S^n\subseteq \mathbb{R}^n\}$ be the set of all possible values of $T$ for different realized samples. Suppose that the underlying problem is to estimate $\lambda$ by maximizing some generalized likelihood function $\mathcal{\mathcal{L_G}}$ (which may be either of \eqref{E8}, \eqref{E9} or \eqref{E10}).
\begin{definition}[Generalized sufficient statistic]\label{GSS}
A statistic $T$ is said to be a sufficient statistic for $\lambda$ with respect to $\mathcal{L_G}$ if $[\mathcal{\mathcal{L_G}}(r_1^n; \lambda)-\mathcal{\mathcal{L_G}}(s_1^n; \lambda) ]$
is independent of $\lambda$ whenever $T(r_1^n)=T(s_1^n)$. $\hfill\Box$
\end{definition}
In the next two propositions, we present the necessary and sufficient condition for a statistic to be sufficient with respect to a generalized likelihood function. Their proofs are available in \cite{Generalized_sufficiency}. 
\begin{proposition}[Generalized factorization theorem]\label{GFT}
A statistic $T$ is said to be a sufficient statistic for $\lambda$ with respect to $\mathcal{L_G}$ if and only if there exist functions $p: \Lambda \times \mathcal{T} \rightarrow \R$ and $q: S^n \rightarrow \R$ such that 
\begin{equation}
    \mathcal{\mathcal{L_G}}(y_1^n; \lambda) = p(\lambda, T(y_1^n)) + q(y_1^n).
\end{equation} 
\end{proposition}
\begin{proposition}\label{GFD}
 Let the deformed probability density associated with the generalized likelihood function $\mathcal{\mathcal{L_G}}$ be given  by
\begin{equation}\label{E_d}
    \widetilde{f}_\lambda(y_1^n) = \frac{\exp[\mathcal{L_G}(y_1^n; \lambda)]}{\int \exp[\mathcal{L_G}(r_1^n; \lambda)]dr_1^n},
\end{equation}
provided the denominator is finite. Then, a statistic $T$ is sufficient for $\lambda$ with respect to $\mathcal{L_G}$ if and only if 
\begin{equation}
\begin{aligned}
  \widetilde{f}_{\lambda_{y_1^n|t}}(y_1^n|t) & := \frac{\widetilde f_\lambda (y_1^n) \textbf{1}(T(y_1^n)=t)}{\widetilde{g}_\lambda(t)}  \\
  & = \begin{cases}
     \frac{\widetilde f_\lambda (y_1^n)}{\widetilde{g}_\lambda(t)} & \text{if $y_1^n\in C_t$}\\
    0 & \text{otherwise},
    \end{cases}
\end{aligned}
\end{equation}
is independent of $\lambda$, where $C_t = \{y_1^n : T(y_1^n)=t\}$ and
$\widetilde{g}_\lambda(t) := \int_{C_t} \widetilde f_\lambda (r_1 ^n)dr_1^n$. $\hfill\Box$
\end{proposition}
We would like to note here that such deformed probability densities $\widetilde{f}_\lambda$, defined in \eqref{E_d}, are also used for generalized (robust) Bayes inference in defining suitable pseudo-posteriors \cite{GB, GMB}. Let us now write $\widetilde{E}_\lambda[\cdot], \widetilde{Var}_\lambda[\cdot]$ and $\widetilde{Cov}_\lambda[\cdot, \cdot]$ to represent, respectively, the mean, the variance and the covariance measures with respect to this deformed density $\widetilde{f}_\lambda$. Next, we give the definition of minimal sufficient statistic in the context of generalized likelihoods. .

\begin{definition}[Minimal sufficient statistic]\label{MSS}
A statistic $T$ is said to be minimal sufficient statistic if $T$ is a function of any other sufficient statistics. Equivalently, $T$ is minimal if, for any two i.i.d. samples $r_1^n$ and $s_1^n$, $T(r_1^n)=T(s_1^n)$ holds if and only if $[\mathcal{L_G}(r_1^n, \lambda) - \mathcal{L_G}(s_1^n, \lambda)]$ is independent of $\lambda$.  $\hfill\Box$ 
\end{definition}
\hspace*{0.02 in}Now we present the formal definition of the exponential family of probability distributions which contains many important examples as special cases (\cite{exp}). This family can be viewed as the projection of KLD on a set  of probability distributions determined by linear constraints, as shown in \cite{Csiszar}. 
\begin{definition}[Exponential family]\label{D2}
Let $ \{ f_\lambda : \lambda \in \Lambda \}$ be a family of probability distributions, where $\Lambda$ is an open subset of $\R^k$. Then, this family is said to form an exponential family, characterized by $w,h,f,\Lambda$ and $S$, if
\begin{equation}\label{E11}
  f_\lambda(y)=\begin{cases}
    \exp\left[h(y)+Z(\lambda)+{w(\lambda)}^T f(y)\right], & \text{for $y \in S$,}\\
    0, & \text{otherwise},
  \end{cases}
\end{equation}
where $S \subseteq \R$ is the (common) support of $f_\lambda$; $e^{Z(\lambda)}$ is the normalizing constant; $w=[w_1, \dots ,w_d]^T$ and $w_i: \Lambda \rightarrow \R$ is differentiable, for $i= 1, \dots ,d$; $f=[f_1, \dots ,f_d]^T$ with $f_i:\R \rightarrow \R$, for $i= 1, \dots ,d$, and $h:\R\rightarrow \R$. We denote this family of distributions by $\mathcal{E} $. $\hfill\Box$  
\end{definition}
Like the exponential family, there are some other generalized families that appeared as the projections of other divergences. These generalized families often contain new important probability distributions along with the distributions that are already there in the exponential family. Let us briefly note the definition and basic properties of two such related families. 
\begin{definition}[$\mathcal{B}^{(\alpha)}$-family]\label{D3}
Let $ \{ f_\lambda : \lambda \in \Lambda \}$ be a family of probability distributions, where $\Lambda$ is an open subset of $\R^k$. Then, this family is said to form a $k$-parameter $ {\mathcal{B}^{(\alpha)}}$-family, characterized by $w,h,f,\Lambda$ and $S$,  if \begin{equation}\label{E12}
  f_\lambda(y)=\begin{cases}
    \left[h(y)+Z(\lambda)+ {w(\lambda)}^T f(y)\right]^\frac{1}{\alpha-1}, & \text{for $y \in S$,}\\
    0, & \text{otherwise},
  \end{cases}
\end{equation}
where $w$, $h$, $f$ and $S$ are the same as in Definition \ref{D2}, and $Z(\lambda) : \Lambda \rightarrow \R$ is the normalizing factor.
\end{definition} 
\begin{definition}[Regular $\mathcal{B}^{(\alpha)}$-family]\label{Dr_b}
An $\mathcal{B}^{(\alpha)}$-family is said to be regular if
\begin{enumerate}
    \item [$(a)$] $S$ does not depend on  $\lambda$;
\item [$(b)$] $d=k$;
\item [$(c)$] $1, w_1, \dots , w_d$  are linearly independent;
\item [$(d)$]  $1,f_1, f_2,\dots,f_d $  are linearly independent.
\end{enumerate}
\end{definition}
\begin{remark}\label{R1}
The following observations can be made regarding the $\mathcal{B^{(\alpha)}}$-family of probability distributions.
\begin{enumerate}
\item [$(a)$] If $f(y)= (1-\alpha)\widetilde{f}(y)$, $Z(\lambda) = (1-\alpha)Z'(\lambda) $ and $h(y)=q^{\alpha-1}(y)$ for some function $\widetilde{f}$, $Z'(\lambda)$, $q$ ($> 0$), and $\alpha \rightarrow 1$, then $\mathcal{B}^{(\alpha)}$-family coincides with $\mathcal{E}$-family;
\item [$(b)$] The $ {\mathcal{B}^{(\alpha)}}$-family can be obtained by the projection of DPD on a set of probability distributions determined by suitable linear constraints (\cite{B-family}). A more general form of this family was considered in Csisz\'{a}r and Mat\'{u}\v{s} \cite{B_family} in connection with the projection problems of Bregman divergences. $\hfill\Box$
\end{enumerate} 
\end{remark}
\begin{definition}[$\mathcal{M}^{(\alpha)}$-family]\label{D4}
Let $ \{ f_\lambda : \lambda \in \Lambda \}$ be a family of probability distributions, where $\Lambda$ is an open subset of $\R^k$. Then, this family is said to form a $k$-parameter $ {\mathcal{M}^{(\alpha)}}$-family, characterized by $w,h,f,\Lambda$ and $S$, if 
\begin{equation}\label{E13}
  f_\lambda(y)=\begin{cases}
    N(\lambda)\left[h(y)+ {w(\lambda)}^T f(y)\right]^\frac{1}{\alpha-1}, & \text{for $y \in S$}\\
    0, & \text{otherwise},
  \end{cases}
\end{equation}
where $w$, $h$, $f$ and $S$ are the same as in Definition \ref{D2}, and  $N(\lambda)^{-1} = \bigintss_S \left[h(y)+ {w(\lambda)}^T f(y)\right]^\frac{1}{\alpha-1}dy$.
\end{definition}

\begin{definition}[Regular $\mathcal{M}^{(\alpha)}$-family]\label{D6}
An $\mathcal{M}^{(\alpha)}$-family is said to be regular if
\begin{enumerate}
    \item [$(a)$] $S$ does not depend on  $\lambda$;
\item [$(b)$] $d=k$;
\item [$(c)$] $1, w_1, \dots , w_d$  are linearly independent;
\item [$(d)$]  $f_1, f_2,\dots,f_d $  are linearly independent.
\end{enumerate}
\end{definition}
\begin{remark}\label{R2} 
The following observations can be made regarding the $\mathcal{M^{(\alpha)}}$-family of probability distributions.
\begin{enumerate}
    \item [$(a)$] If $f(y)= (1-\alpha)\widetilde{f}(y)$ and $h(y)=q^{\alpha-1}(y)$ for some function $\widetilde{f}$, $q$ ($> 0$) and $\alpha \rightarrow 1$, then $\mathcal{M}^{(\alpha)}$-family coincides with $\mathcal{E}$-family;
    \item [$(b)$] The $ {\mathcal{M}^{(\alpha)}}$-family can be obtained by the projection of LDPD on a set of probability distributions determined by suitable linear constraints (\cite{Forward_projection, Reverse_projection}); 
    \item[$(c)$] Any member of $ {\mathcal{B}^{(\alpha)}}$-family  can be expressed as a member of $ {\mathcal{M}^{(\alpha)}}$-family and vice versa. However, this is not true for regular $ {\mathcal{B}^{(\alpha)}}$ and $ {\mathcal{M}^{(\alpha)}}$-families (\cite{Projection_theorems}).~$\hfill\Box$
\end{enumerate} 
\end{remark}

Finally, let us formally define the family of Student's distributions from \cite{JV}.
\begin{definition}[Student's distribution]\label{s_d}
A family of Student's distributions is the set of probability distributions $\mathcal{P} =\{f_{\mu,\sigma^2} : (\mu, \sigma^2) \in \R \times \R^+\} $ with the pdf of the form
$$f_{\mu,\sigma^2}(y)= F_\nu(\sigma^2)\left[1+\frac{(y-\mu)^2}{\nu \sigma^2}\right]_+^{-\frac{\nu+1}{2}},$$
where 
\begin{eqnarray*}
 F_\nu(\sigma^2)   = \begin{cases}
        \frac{\Gamma(\frac{\nu+1}{2})}{\Gamma(\frac{\nu}{2})\sqrt{\pi \nu \sigma^2}} & \text{if} ~ \nu \in (0,\infty)\\
        \frac{\Gamma\left(1-\frac{\nu}{2}\right)}{\Gamma\left(1-\frac{\nu+1}{2}\sqrt{-\pi \nu \sigma^2}\right)} & \text{if}~ \nu \in (-\infty,0),
    \end{cases}
\end{eqnarray*} 
and the support is given by
\begin{eqnarray*}
    S = \begin{cases} 
          \R, & \text{if} ~ \nu \in (0,\infty)\\
          \left\{y : \left(\frac{y-\mu}{\sigma}\right)^2 < -\nu \right\} & \text{if} ~ \nu \in (-\infty, 0),
\end{cases}
\end{eqnarray*}
and,
 for any $r\in \mathbb{R}$, $[r]_+ :=\max\{r,0\}$. Moreover, this family of probability distributions is commonly referred to as the family of Student's $t$-distributions (resp. the family of Student's $r$-distributions) when $\nu>0$ (resp. $\nu<0$).
\end{definition}
\section{Estimating equation and generalized likelihood function associated with the minimum LNRE estimation}\label{S3}
In this section, we first consider the minimization of the LNRE with respect to the parameter $\lambda$ and obtain the associated estimating equations and the corresponding likelihood function, which is further used in subsequent sections in the study of generalized sufficiency.

Let us again consider the set-up and notation of parametric estimation problem, as mentioned in Section \ref{S1}, for the family $\mathcal{P}$ with the common support $S \subseteq \mathbb{R}$. Further, our interest lies in estimating the parameter $\lambda$ for which $g$ and $f_\lambda$ are closest, measuring this closeness in terms of the LNRE, $\mathcal{RE}_{\alpha,\beta}^{\mathcal{LN}}(g, f_\lambda)$, with respect to $\lambda$, i.e., 
\begin{equation*}
min_{\lambda \in \Lambda} \mathcal{RE}_{\alpha,\beta}^{\mathcal{LN}} (g,f_\lambda)  = min_{\lambda \in \Lambda} \frac{\alpha}{ \beta(\beta-\alpha)} \log \int g^\beta(y) f_\lambda^{\alpha-\beta}(y)dy - \frac{1}{\beta-\alpha} \log \int g^\alpha(y)dy  + \frac{1}{\beta}\log \int f_\lambda^\alpha(y)dy.
\end{equation*}
By routine differentiation, the corresponding estimating equation is given by
\begin{equation*}
\frac{\int {g}^\beta(y) f_\lambda^{\alpha-\beta}(y)u(y,\lambda)dy}{\int {g}^\beta(y) f_\lambda^{\alpha-\beta}(y)dy}  = \frac{\int f_\lambda^\alpha(y)u(y,\lambda)dy}{\int f_\lambda^\alpha(y)dy}.
\end{equation*}
But, in practice, the true distribution of the data is not known. So, we consider $\widehat{g}$, a suitable estimate of $g$ based on the data. The above estimating equation can then be written in the sample form as
\begin{equation}\label{E14}
\frac{\frac{1}{n} \sum_{j=1} ^ {n} \widehat{g}^{\beta-1}(y_j) f_\lambda^{\alpha-\beta}(y_j)u(y_j,\lambda)}{\frac{1}{n} \sum_{j=1} ^ {n} \widehat{g}^{\beta-1}(y_j) f_\lambda^{\alpha-\beta}(y_j)} = \frac{\int f_\lambda^\alpha(y)u(y,\lambda)dy}{\int f_\lambda^\alpha(y)dy}.
\end{equation}
Note that if the underlying distribution is discrete, then the empirical measure of the sample can be used to construct $\widehat{g}$. However, for continuous distribution, $\widehat{g}$ will be some nonparametric (e.g., kernel type) estimate based on sample.
We call the estimator obtained by solving the above (sample) estimating equation to be a MLNREE. This estimating equation is a two-parameter generalization of the estimating equation associated with the LDPD, i.e., the expression given in \eqref{E7} coincides with the expression given in \eqref{E14} when $\beta =1$. 
Further, for $\alpha>0 $ and $\beta \in \R$, we consider the generalized likelihood function
\begin{eqnarray}
 \mathcal{L}_{\mathcal{RE}}^{(\alpha,\beta)}(y_1^n;\lambda) &:=& \frac{1}{\alpha-\beta}\log \frac{1}{n}\sum_{j=1} ^ {n} \widehat{g}^{\beta-1}(y_j) f_\lambda^{\alpha-\beta}(y_j) - \frac{1}{\alpha}\log \int f_\lambda^\alpha(y)  dy \label{E16} \\ 
 &=& \frac{1}{\alpha-\beta}\log \frac{1}{n}\sum_{j=1} ^ {n} \widehat{g}^{\beta-1}(y_j) f_\lambda^{\alpha-\beta}(y_j) - \log || f_\lambda ||_\alpha \label{Ex} \\ 
 &=& \log \frac{ \left(\frac{1}{n}\sum_{j=1} ^ {n}\widehat{g}^{\beta-1}(y_j) f_\lambda^{\alpha-\beta}(y_j)\right)^{\frac{1}{\alpha-\beta}}}{ || f_\lambda ||_\alpha }\label{Exx},
\end{eqnarray}
where $|| f_\lambda ||_\alpha = \left[\int f_\lambda^\alpha(y)  dy\right]^\frac{1}{\alpha}$. Note that the MLNREE, obtained by solving \eqref{E14}, can also be derived by maximizing $\mathcal{L}_{\mathcal{RE}}^{(\alpha,\beta)}(y_1^n;\lambda)$ with respect to $\lambda$, since the first order optimality condition of this maximization problem is the same as in \eqref{E14}. In particular, it coincides with the usual log-likelihood $\mathcal{L}(y_1^n,\lambda)$ in \eqref{E8} as $\alpha\rightarrow 1$ and $\beta \rightarrow 1$. Thus, $\mathcal{L}_{\mathcal{RE}}^{(\alpha,\beta)}(y_1^n;\lambda)$ can be indeed considered as a generalized likelihood function associated with the LNRE. Now, based on this generalized likelihood function, $\mathcal{L}_{\mathcal{RE}}^{(\alpha,\beta)}(y_1^n;\lambda)$, we can obtain the generalized sufficient statistic for a given family of probability distributions. However, the notion of sufficient statistic is useful when it is fixed in number (independent of sample size). Thus, in the next section, we study a family of probability distributions that has a fixed number of sufficient statistics.

\section{Generalized Fisher-Darmois-Koopman-Pitman theorem}\label{S4}
 Similar to the exponential family, we can obtain a two-parameter power-law family of distributions as a projection of LNRE on a set  of probability distributions determined by suitable linear constraints (\cite{LNRE}, \cite{Scale_invariant}). We call this family as the $\mathcal{M}^{(\alpha, \beta)}$-family, which is formally defined below.
\begin{definition}\label{D5}
Let $ \{ f_\lambda : \lambda \in \Lambda \}$ be a family of probability distributions, where $\Lambda$ is an open subset of $\R^k$. Then, this family is said to form a $k$-parameter $ {\mathcal{M}^{(\alpha,\beta)}}$-family, characterized by $w,h,f,\Lambda$ and $S$, if 
\begin{equation}\label{E17}
  f_\lambda(y)=\begin{cases}
    N(\lambda){[h(y)+{w(\lambda)}^T f(y)]}^\frac{1}{\alpha-\beta}, & \text{for $y\in S$}.\\
    0 & \text{otherwise},
  \end{cases}
\end{equation}
where $w$, $h$, $f$ and $S$ are the same as in Definition \ref{D2}, and  ${N(\lambda)}^{-1} = \bigints_S {\left[h(x)+{w(\lambda)}^T f(x)\right]}^\frac{1}{\alpha-\beta} dx$.
\end{definition} 
\begin{remark}\label{R3}
The following observations can be made regarding the $\mathcal{M^{(\alpha,\beta)}}$-family of probability distributions.
\begin{enumerate}
\item[$(a)$] $\mathcal{E} \subset \mathcal{M}^{(\alpha)} \subset \mathcal{M}^{(\alpha, \beta)} $;
\item[$(b)$] If $f(y)= (\alpha-\beta)\widetilde{f}(y)$ and $h(y)=q^{\alpha-\beta}(y)$ for some functions $\widetilde{f}$ and $ q $ ($> 0$) and $\alpha \rightarrow \beta$, then $\mathcal{M}^{(\alpha, \beta)}$- family coincides with $\mathcal{E}$-family;

\item[$(c)$] For $\beta=1$, $ {\mathcal{M}^{(\alpha,\beta)}}$-family coincides with $\mathcal{M^{(\alpha)}}$-family;

\item[$(d)$] The $\mathcal{M}^{(\alpha,\beta)}$-family contains both Student's $t$-distribution and Student's $r$-distribution as special cases (see Example \ref{Ex1} and Section \ref{S7.2} later) . $\hfill\Box$
\end{enumerate}
\end{remark}
\hspace*{0.02 in}Below we give the definition of regular ${\mathcal{M}^{(\alpha,\beta)}}$-family in the line of Definition \ref{D6}.
\begin{definition}
An $\mathcal{M}^{(\alpha, \beta)}$-family is said to be regular if
\begin{enumerate}
    \item [$(a)$] $S$ does not depend on  $\lambda$;
\item [$(b)$] $d=k$;
\item [$(c)$] $1, w_1, \dots , w_d$  are linearly independent;
\item [$(d)$]  $f_1, f_2,\dots,f_d $  are linearly independent. $\hfill\Box$
\end{enumerate}
\end{definition}
In the following theorem, we show that, if a given family has a fixed number of sufficient statistics, independent of sample size, with respect to the likelihood function in \eqref{E16}-\eqref{Exx} associated with the LNRE, then the family is indeed  the ${\mathcal{M}^{(\alpha,\beta)}}$-family.  
\begin{theorem}\label{T1}
 Let $\mathcal{P}= \{f_\lambda: \lambda\in\Lambda\subset \R^k\}$ be a family of probability distributions with common support $S$, an open subset of $\R$. Let $y_1^n$ be an i.i.d. sample drawn from some $f_\lambda \in \mathcal{P} $. For some $t\;(k \leq t < n)$, let $T_i (y_1^n)$ be differentiable functions on $S^n$, for $i= 1, \dots ,t$. If ${[T_1,\dots,T_t]}^T$ is a sufficient statistic for $\lambda$ with respect to the likelihood function $\mathcal{L}_{\mathcal{RE}}^{(\alpha,\beta)}$ given in \eqref{Ex}, then there exist functions $h:S\rightarrow \R, w_i: \lambda \rightarrow \R$, and $f_i:S \rightarrow \R$, for $i= 1, \dots ,d$, such that $f_\lambda$ can be expressed as \eqref{E17}.
\end{theorem}
\noindent\textbf{Proof :} We consider the following two cases.
\\\textbf{Case-I:}
Let $k=t=1$. Taking derivative of the likelihood function given in \eqref{Ex} with respect to  $\lambda$ on both sides, we get
\begin{equation}\label{E18}
\frac{\partial \mathcal{L}_{\mathcal{RE}}^{(\alpha,\beta)}(y_1^n;\lambda) }{\partial \lambda} =  \frac{1}{\alpha-\beta} \frac{\frac{1}{n}\sum_{j=1}^{n} \widehat{g}^{\beta-1}(y_j) \frac{\partial{f_\lambda ^{\alpha-\beta}(y_j) }}{\partial \lambda}}{\frac{1}{n}\sum_{j=1}^{n} \widehat{g}^{\beta-1}(y_j) f_\lambda ^{\alpha-\beta} (y_j)}-
\frac{\partial{\log||f_\lambda||_\alpha}}{\partial{\lambda}}. 
\end{equation}

\noindent If $T$ is a sufficient statistic for $\lambda$ with respect to $\mathcal{L}_{\mathcal{RE}}^{(\alpha,\beta)}(y_1^n;\lambda)$, then, by Proposition~\ref{GFT}, we have
\begin{equation}\label{E19}
    \mathcal{L}_{\mathcal{RE}}^{(\alpha,\beta)}(y_1^n;\lambda) = p(\lambda, T)+q(y_1^n),
\end{equation}
 for some functions $p$ and $q$, where $q$ is independent of $\lambda$. Consequently, we have
 \begin{equation}\label{E20}
     \frac{\partial \mathcal{L}_{\mathcal{RE}}^{(\alpha,\beta)}(y_1^n;\lambda) }{\partial \lambda} = \frac{\partial [p(\lambda,T)] }{\partial \lambda} = p_1(\lambda, T) \quad \text(say).
 \end{equation}
 By comparing \eqref{E18} and \eqref{E20}, we get
 \begin{equation}\label{E21}
 p_1(\lambda, T)= \frac{1}{\alpha-\beta} \frac{\frac{1}{n}\sum_{j=1}^{n} \widehat{g}^{\beta-1}(y_j) \frac{\partial{f_\lambda ^{\alpha-\beta}(y_j) }}{\partial \lambda}}{\frac{1}{n}\sum_{j=1}^{n} \widehat{g}^{\beta-1}(y_j) f_\lambda ^{\alpha-\beta} (y_j)}-
\frac{\partial{\log||f_\lambda||_\alpha}}{\partial{\lambda}}. 
 \end{equation}
Now, for a given particular value of $\lambda \in \Lambda$, $\lambda=\lambda'$ (say), the above equation becomes 
 \begin{equation}\label{E22}
 \begin{aligned}
 p_1(\lambda', T) & = \frac{1}{\alpha-\beta} \frac{\sum_{j=1}^{n} \widehat{g}^{\beta-1}(y_j) \left[\frac{\partial{f_\lambda ^{\alpha-\beta}(y_j) }}{\partial \lambda}\right]_{\lambda'}}{\sum_{j=1}^{n} \widehat{g}^{\beta-1}(y_j) f_{\lambda'} ^{\alpha-\beta} (y_j)}- 
\left[\frac{\partial{\log||f_\lambda||}}{\partial{\lambda}}\right]_{\lambda'} \\
& =\frac{1}{\alpha-\beta} \frac{\sum_{j=1}^{n} g(y_j)}{\sum_{j=1}^{n} h(y_j)}- 
\left[\frac{\partial{\log||f_\lambda||_\alpha}}{\partial{\lambda}}\right]_{\lambda'},
\end{aligned}
 \end{equation}
where $g(y_j) = \widehat{g}^{\beta-1}(y_j) \left[\frac{\partial{f_\lambda ^{\alpha-\beta}(y_j) }}{\partial \lambda}\right]_{\lambda'}$ and $h(y_j)= \widehat{g}^{\beta-1}(y_j) f_{\lambda'} ^{\alpha-\beta} (y_j)$. Further, let  $\widehat{T}_1 =\sum_{j=1}^{n} g(y_j) $ and $\widehat{T}_2 =\sum_{j=1}^{n} h(y_j) $. 
Then, in view of \eqref{E22}, there exists a function $\xi$ such that
\begin{equation}\label{E23}
    T=\xi\left(\widehat{T}_1,\widehat{T}_2\right).
\end{equation} 
Consequently, \eqref{E20} can be rewritten as
\begin{equation}\label{E24}
  \frac{\partial}{\partial \lambda}\mathcal{L}_{\mathcal{RE}}^{(\alpha,\beta)}(y_1^n;\lambda)  = p_2\left(\lambda,\widehat{T}_1,\widehat{T}_2\right),  
\end{equation}
for some function $p_2$. Further, from  \eqref{Exx}, we have
$$ \exp\left[(\alpha-\beta) \mathcal{L}_{\mathcal{RE}}^{(\alpha,\beta)}(y_1^n;\lambda) \right]= \frac{1}{n}\sum_{j=1}^{n}\frac{ \widehat{g}^{\beta-1}(y_j) f_\lambda ^{\alpha-\beta} (y_j)}{ 
{||f_\lambda||}^{\alpha-\beta}},$$
which, in view of \eqref{E19}, can equivalently be written as
\begin{equation}\label{E26}
  \frac{1}{n}\sum_{j=1}^{n}\frac{ \widehat{g}^{\beta-1}(y_j) f_\lambda ^{\alpha-\beta} (y_j)}{ 
{||f_\lambda||}^{\alpha-\beta}} = \exp\left[(\alpha-\beta) p(\lambda, T)\right]\exp\left[(\alpha-\beta) q(y_1^n)\right].
\end{equation}
Again, for $\lambda=\lambda'$, the above equation becomes 
\begin{equation*}
\widehat{T}_2= n{||f_{\lambda'}||}^{\alpha-\beta}\exp\left[(\alpha-\beta) p(\lambda', T)\right]\exp\left[(\alpha-\beta)q(y_1^n)\right].
\end{equation*}
Since $T=\xi\left(\widehat{T}_1,\widehat{T}_2\right)$, we conclude from the above equality that $q(y_1^n)$ must be a function of $\widehat{T}_1$ and $\widehat{T}_2$ only. Consequently, \eqref{E26} can be written as
\begin{equation}\label{E27}
\sum_{j=1}^{n}\widehat{g}^{\beta-1}(y_j) f_\lambda ^{\alpha-\beta} (y_j)= p_3\left(\lambda,\widehat{T}_1,\widehat{T}_2\right),
\end{equation}
where $$p_3\left(\lambda,\widehat{T}_1,\widehat{T}_2\right)= n {||f_{\lambda}||}^{\alpha-\beta}\exp[(\alpha-\beta) p(\lambda, T)]\exp[(\alpha-\beta)q(y_1^n)].$$ Taking partial derivative with respect to $y_j$ on both sides of  the expression given in \eqref{E27}, we get
\begin{equation}\label{E28}
    \begin{aligned}
        \frac{\partial}{\partial y_j}\left[\widehat{g}^{\beta-1}(y_j) f_\lambda ^{\alpha-\beta} (y_j)\right] & = \frac{\partial}{\partial \widehat{T}_1} p_3\left(\lambda,\widehat{T}_1,\widehat{T}_2\right) \frac{\partial \widehat{T}_1}{\partial y_j} + \frac{\partial}{\partial \widehat{T}_2} p_3\left(\lambda,\widehat{T}_1,\widehat{T}_2\right) \frac{\partial \widehat{T}_2}{\partial y_j}   \\
& = \frac{\partial}{\partial \widehat{T}_1} p_3\left(\lambda,\widehat{T}_1,\widehat{T}_2\right) \frac{\partial {g(y_j)}}{\partial y_j} + \frac{\partial}{\partial \widehat{T}_2} p_3\left(\lambda,\widehat{T}_1,\widehat{T}_2\right) \frac{\partial {h(y_j)}}{\partial y_j} \\
& = H(y_j,\lambda)\quad \text{(say)}.
    \end{aligned}
\end{equation}
Again, taking partial derivative with respect to $y_i$, for $i \neq j $, on both sides of the above equation, we get
\begin{equation}\label{E29}
  \frac{\partial}{\partial y_i}\left[\frac{\partial}{\partial \widehat{T}_1} p_3\left(\lambda,\widehat{T}_1,\widehat{T}_2\right)\right] \frac{\partial {g(y_j)}}{\partial y_j} + \frac{\partial}{\partial y_i}\left[\frac{\partial}{\partial \widehat{T}_2} p_3\left(\lambda,\widehat{T}_1,\widehat{T}_2\right)\right] \frac{\partial {h(y_j)}}{\partial y_j} = 0.
\end{equation}
 Since both $\frac{\partial {g(y_j)}}{\partial y_j}$ and $\frac{\partial {h(y_j)}}{\partial y_j}$ are non zero, we can say from the above equation that 
\begin{equation}\label{E30}
\frac{\partial}{\partial y_i}\left[\frac{\partial}{\partial \widehat{T}_1} p_3\left(\lambda,\widehat{T}_1,\widehat{T}_2\right)\right] = 0 \quad \text{if and only if} \quad
\frac{\partial}{\partial y_i}\left[\frac{\partial}{\partial \widehat{T}_2} p_3\left(\lambda,\widehat{T}_1,\widehat{T}_2\right)\right] = 0,
\end{equation} 
which can equivalently be written as
\begin{equation*}
\frac{\partial}{\partial y_i}\left[\frac{\partial}{\partial \widehat{T}_1} p_3\left(\lambda,\widehat{T}_1,\widehat{T}_2\right)\right] \neq 0 \quad \text{if and only if} \quad
\frac{\partial}{\partial y_i}\left[\frac{\partial}{\partial \widehat{T}_2} p_3\left(\lambda,\widehat{T}_1,\widehat{T}_2\right)\right] \neq 0.
\end{equation*} 
 Let $\frac{\partial}{\partial y_i}\left[\frac{\partial}{\partial \widehat{T}_1} p_3\left(\lambda,\widehat{T}_1,\widehat{T}_2\right)\right] \neq 0$. Then, from \eqref{E29}, we get
 \begin{equation}\label{E31}
 \begin{aligned}
       \frac{\frac{\partial}{\partial y_i}\left[\frac{\partial}{\partial \widehat{T}_1} p_3\left(\lambda,\widehat{T}_1,\widehat{T}_2\right)\right] }{\frac{\partial}{\partial y_i}\left[\frac{\partial}{\partial \widehat{T}_2} p_3\left(\lambda,\widehat{T}_1,\widehat{T}_2\right)\right]} & = - \frac{\frac{\partial {h(y_j)}}{\partial y_j}}{\frac{\partial {g(y_j)}}{\partial y_j}} = \phi(y_j) \quad \text{(say)}. \\
 \end{aligned}  
 \end{equation}  
Note that the right-hand side of above equation is a function of $y_j$ only. Thus, both $\frac{\partial}{\partial \widehat{T}_1} p_3\left(\lambda,\widehat{T}_1,\widehat{T}_2\right)$ and $\frac{\partial}{\partial \widehat{T}_2} p_3\left(\lambda,\widehat{T}_1,\widehat{T}_2\right)$ must be in the form of

\begin{equation*}
\frac{\partial}{\partial \widehat{T}_1} p_3\left(\lambda,\widehat{T}_1,\widehat{T}_2\right)= \zeta(y_1^n,\lambda)A_1(y_j) + B_1(y_1, \dots ,y_{i-1},y_{i+1},\dots,y_n,\lambda)
\end{equation*}
and 
\begin{equation*}
\frac{\partial}{\partial \widehat{T}_2} p_3\left(\lambda,\widehat{T}_1,\widehat{T}_2\right)=
\zeta(y_1^n,\lambda)A_2(y_j) + B_2(y_1, \dots ,y_{i-1},y_{i+1},\dots,y_n,\lambda),
\end{equation*}
for some real valued functions $\zeta, A_1, A_2, B_1 $ and $B_2$. Upon using this in \eqref{E28}, we get
\begin{equation}\label{E32}
H(y_j,\lambda) = \zeta(y_1^n,\lambda) H_1(y_j) + H_2(y_j,\lambda),
\end{equation}
where
$H_1(y_j) = A_1(y_j)\frac{\partial}{\partial y_j} [g(y_j)] + A_2(y_j)\frac{\partial}{\partial y_j} [h(y_j)] $ and $ H_2(y_j) = B_1(y_j,\lambda)\frac{\partial}{\partial y_j} [g(y_j)] + B_2(y_j)\frac{\partial}{\partial y_j} [h(y_j)]  $.
Note that both sides of \eqref{E32} agree if $\zeta$ is a function of $y_j$ and $\lambda$ only. Further, the equation given in \eqref{E31} holds for all $i \neq j$. This implies that both $B_1$ and $B_2$ must be function of $y_j$ and $\lambda$ only. Consequently,
\begin{equation}\label{E33}
 \frac{\partial}{\partial \widehat{T}_1} p_3\left(\lambda,\widehat{T}_1,\widehat{T}_2\right) \hspace{0.2 in} and \hspace{0.2 in}  \frac{\partial}{\partial \widehat{T}_2} p_3\left(\lambda,\widehat{T}_1,\widehat{T}_2\right)\end{equation}
are functions of $y_j$ and $\lambda$ only. Further, this implies that $\frac{\partial}{\partial y_i}\left[\frac{\partial}{\partial \widehat{T}_1} p_3\left(\lambda,\widehat{T}_1,\widehat{T}_2\right)\right] = 0$ and $\frac{\partial}{\partial y_i}\left[\frac{\partial}{\partial \widehat{T}_2} p_3\left(\lambda,\widehat{T}_1,\widehat{T}_2\right)\right] = 0$, for all $i \neq j$. Again, these imply that both $\frac{\partial}{\partial \widehat{T}_1} p_3\left(\lambda,\widehat{T}_1,\widehat{T}_2\right)$ and $\frac{\partial}{\partial \widehat{T}_2} p_3\left(\lambda,\widehat{T}_1,\widehat{T}_2\right)$ are independent of $y_i\text{'}s$ for $i \neq j$. Since $p_3\left(\lambda,\widehat{T}_1,\widehat{T}_2\right)$ is a function of $\widehat{T}_1, \widehat{T}_2$ and $\lambda$ only and both $\widehat{T}_1 $ and $ \widehat{T}_2$ are symmetric functions of $y_j\text{'}s$, $j=1, \dots ,n$, we must have that both elements in \eqref{E33} are function of $\lambda$ only. Suppose \begin{equation*}
 \frac{\partial}{\partial \widehat{T}_1} p_3\left(\lambda,\widehat{T}_1,\widehat{T}_2\right) = u_1(\lambda) \hspace{0.2 in} and \hspace{0.2 in}  \frac{\partial}{\partial \widehat{T}_2} p_3\left(\lambda,\widehat{T}_1,\widehat{T}_2\right) = u_2(\lambda),
 \end{equation*}
for some functions $u_1$ and $u_2$ of $\lambda$. Then, we have 
\begin{equation}\label{E34}
p_3\left(\lambda,\widehat{T}_1,\widehat{T}_2\right) = u_1(\lambda) \widehat{T}_1+
u_2(\lambda)\widehat{T}_2+u_3(\lambda),  
 \end{equation}
for some $u_3(\lambda)$. Now, by using \eqref{E27} in the above equation, we get
$$ \sum_{j=1}^{n} \widehat{g}^{\beta-1}(y_j) f_\lambda ^{\alpha-\beta} (y_j) = u_1(\lambda) \sum_{j=1}^{n} g(y_j)+
u_2(\lambda) \sum_{j=1}^{n} h(y_j)+u_3(\lambda).  $$
Further, this gives
\begin{equation*}
 f_\lambda(y_j) = [u_1(\lambda) g^*(y_j)+
u_2(\lambda)h^*(y_j)+u_3(\lambda)v^*(y_j)]^\frac{1}{\alpha-\beta} ,
\end{equation*}
where $ g^*(y_j) =  \left[\frac{\partial{f_\lambda ^{\alpha-\beta}(y_j) }}{\partial \lambda}\right]_{\lambda'}$, $ h^*(y_j) =  f_{\lambda'} ^{\alpha-\beta} (y_j)$ and $v^*(y_j) = \frac{1}{n\widehat{g}^{\beta-1}(y_j)}$.
Hence $f_\lambda \in \mathcal{M}^{(\alpha,\beta)}$.\\
\textbf{Case-II:} Let $k<t$. In particular, it suffices to show the result for $k=1$ and $t=2$. Let $(T_1,T_2)$ be a pair of sufficient statistics of $\lambda$ with respect to the likelihood function $\mathcal{L}_{\mathcal{RE}}^{(\alpha,\beta)}$. Then, following similar steps as in Case-I, we get
\begin{equation}
 p_1(\lambda', T_1, T_2) =\frac{1}{\alpha-\beta} \frac{\sum_{j=1}^{n} g_1(y_j)}{\sum_{j=1}^{n} h_1(y_j)}- 
\left[\frac{\partial{\log||f_\lambda||_\alpha}}{\partial{\lambda}}\right]_{\lambda'},
\end{equation}
for particular value of $\lambda = \lambda'$, where $g_1(y_j) = \widehat{g}^{\beta-1}(y_j) \left[\frac{\partial{f_\lambda ^{\alpha-\beta}(y_j) }}{\partial \lambda}\right]_{\lambda'}$ and $h_1(y_j)= \widehat{g}^{\beta-1}(y_j) f_{\lambda'} ^{\alpha-\beta} (y_j)$. Again, for $\lambda = \lambda''$ ($\neq \lambda'$), we have 
\begin{equation}
 p_1(\lambda'', T_1, T_2) =\frac{1}{\alpha-\beta} \frac{\sum_{j=1}^{n} g_2(y_j)}{\sum_{j=1}^{n} h_2(y_j)}- 
\left[\frac{\partial{\log||f_\lambda||_\alpha}}{\partial{\lambda}}\right]_{\lambda''},
\end{equation}
 where $g_2(y_j) = \widehat{g}^{\beta-1}(y_j) \left[\frac{\partial{f_\lambda ^{\alpha-\beta}(y_j) }}{\partial \lambda}\right]_{\lambda''}$ and $h_2(y_j)= \widehat{g}^{\beta-1}(y_j) f_{\lambda''} ^{\alpha-\beta} (y_j)$. Let $\widehat{T}_1= \sum_{j=1} ^ {n} g_1(y_j)$, $\widehat{T}_2= \sum_{j=1} ^ {n} h_1(y_j)$, $\widehat{T}_3= \sum_{j=1} ^ {n} g_2(y_j)$ and $\widehat{T}_4= \sum_{j=1} ^ {n} h_2(y_j)$. Then, we have $p_1(\lambda, T_1, T_2) = \xi \left(\lambda, \widehat{T}_1,\widehat{T}_2,\widehat{T}_3,\widehat{T}_4\right)$, for some $\xi$, where $\widehat{T}_1, \widehat{T}_2, \widehat{T}_3, \widehat{T}_4$ are symmetric functions of $y_j\text{'}s$. By proceeding in the same lines as in Case-I, we get that all $ \frac{\partial}{\partial \widehat{T}_i} p_3\left(\lambda,\widehat{T}_1,\widehat{T}_2,\widehat{T}_3, \widehat{T}_4\right)$, $i= 1,2,3,4$, are function of $\lambda$ only. Further, this gives

\begin{equation*}
    p_3\left(\lambda,\widehat{T}_1,\widehat{T}_2, \widehat{T}_3, \widehat{T}_4\right) = u_1(\lambda) \widehat{T}_1+
u_2(\lambda)\widehat{T}_2+u_3(\lambda) \widehat{T}_3+u_4(\lambda) \widehat{T}_4+u_5(\lambda) ,
\end{equation*}
for some functions $u_1,u_2,u_3,u_4$ and $u_5$. Further, this 
implies that
\begin{equation*}
  \widehat{g}^{\beta-1}(y_j) f_\lambda ^{\alpha-\beta} (y_j) = u_1(\lambda) g_1(y_j)+
u_2(\lambda)h_1(y_j)+u_3(\lambda)g_2(y_j)+u_4(\lambda)h_2(y_j)+\frac{u_5(\lambda)}{n},    
\end{equation*}
or equivalently,
\begin{equation*}
 f_\lambda (y_j) = [u_1(\lambda) g_1^*(y_j)+
u_2(\lambda)h_1^*(y_j)+u_3(\lambda)g_2^*(y_j)+u_4(\lambda)h_2^*(y_j)+u_5(\lambda)v^*(y_j)]^{\frac{1}{\alpha - \beta}},
\end{equation*}
where $g_1^*(y_j) = \left[\frac{\partial{f_\lambda ^{\alpha-\beta}(y_j) }}{\partial \lambda}\right]_{\lambda'},$ $ h_1^*(y_j) = f_{\lambda'} ^{\alpha-\beta} (y_j),$ $ g_2^*(y_j) = \left[\frac{\partial{f_\lambda ^{\alpha-\beta}(y_j) }}{\partial \lambda}\right]_{\lambda''},$ $ h_2^*(y_j) = f_{\lambda''} ^{\alpha-\beta} (y_j)$ and $v^*(y_j)=\frac{1}{n\widehat{g}^{\beta-1}(y_j)} $. Thus, $f_\lambda \in \mathcal{M}^{(\alpha,\beta)}$. Hence the result is proved.$\hfill\Box$
\par
\bigskip
 One natural question arises whether the converse of the above theorem is true or not. In the following theorem, we show that an $\mathcal{M}^{(\alpha,\beta)}$- family always has
a fixed number of sufficient statistics, independent of sample size, with respect to the likelihood function associated with the LNRE.
\begin{theorem}\label{T2}
 Let $\mathcal{P}= \{f_\lambda: \lambda\in\Lambda\subset \R^k\}$ be a family of probability distributions with common support $S$, an open subset of $\R$.
Let $T = [T_1,\dots,T_t]^T$ be $t$ statistics defined on $S^n$, where $t \geq k$. Then, $T$ is a sufficient statistic for $\lambda$ with respect to $\mathcal{L}_{\mathcal{RE}}^{(\alpha,\beta)}$ if 
\begin{enumerate}
    \item [$(i)$] $\mathcal{P}$ is a $k$- parameter $\mathcal{M}^{(\alpha,\beta)}$- family as in \eqref{E17};
    \item [$(ii)$] there exist real-valued functions $\varphi_1, \dots , \varphi_d$, defined on $\mathcal{T}$, such that for any i.i.d. sample $y_1^n$,
\begin{equation}\label{E35}
    \varphi_i (T(y_1^n)) = \frac{\overline{f_i}}{\overline{h}}, \quad i=1,\dots,d,
\end{equation}
 where $\overline{f_i}(y_1^n) = \frac{1}{n}\sum_{j=1}^{n} \widehat{g}^{\beta-1}(y_j) f_i(y_j)$ and $\overline{h}(y_1^n) = \frac{1}{n}\sum_{j=1}^{n} \widehat{g}^{\beta-1}(y_j) h(y_j)$.
\end{enumerate}
\end{theorem}
\noindent\textbf{Proof :} Let $\mathcal{P} = M^{(\alpha,\beta)}$ as in \eqref{E17}. Further, let $r_1 ^n$ and $s_1^n$ be two i.i.d. samples from some $f_\lambda \in M^{(\alpha,\beta)}$ such that $T(r_1 ^n)= T(s_1^n)$. For $i= 1,\dots,d$, let
\begin{equation*}
  \varphi_i (T(r_1^n)) = \frac{\overline{f_i} (r_1 ^n)}{\overline{h}(r_1 ^n)} \hspace{0.2 in}and\hspace{0.2in }
 \varphi_i (T(s_1^n)) = \frac{\overline{f_i} (s_1^n)}{\overline{h}(s_1^n)}.  
\end{equation*}

\noindent Consider
\begin{equation*}
\begin{aligned}
\mathcal{L}_{\mathcal{RE}}^{(\alpha,\beta)}(r_1 ^n;\lambda) - \mathcal{L}_{\mathcal{RE}}^{(\alpha,\beta)}(s_1^n;\lambda)
& = \frac{1}{\alpha-\beta} \log \frac{\frac{1}{n}\sum_{j=1}^{n} \widehat{g}^{\beta-1}(r_j) f_\lambda ^{\alpha-\beta} (r_j)} { \frac{1}{n}\sum_{j=1}^{n} \widehat{g}^{\beta-1}(s_j) f_\lambda ^{\alpha-\beta} (s_j)}.
\end{aligned}
\end{equation*}
After incorporating the expression of $f_\lambda$ in the above equation, we get
\begin{equation}\label{E36}
\begin{aligned}
\mathcal{L}_{\mathcal{RE}}^{(\alpha,\beta)}(r_1 ^n;\lambda) - \mathcal{L}_{\mathcal{RE}}^{(\alpha,\beta)}(s_1^n;\lambda) & =  \frac{1}{\alpha-\beta} \log \frac{\frac{1}{n} \sum_{j=1}^{n} \widehat{g}^{\beta-1}(r_j) {N(\lambda)}^{\alpha-\beta}{\left[h(r_j)+{w(\lambda)}^T f(r_j)\right]}} {\frac{1}{n} \sum_{j=1}^{n} \widehat{g}^{\beta-1}(s_j){N(\lambda)}^{\alpha-\beta}{\left[h(s_j)+{w(\lambda)}^T f(s_j)\right]}} \\
& = \frac{1}{\alpha-\beta} \left[\log \frac{\overline{h}(r_1^n)}{\overline{h}(s_1^n)}
+\log \frac{1+w(\lambda)^T\overline{f}(r_1^n)/\overline{h}(r_1^n)}{ 1+w(\lambda)^T\overline{f}(s_1^n)/\overline{h}(s_1^n)} \right]\\
& = \frac{1}{\alpha-\beta} \log \left[\overline{h}(r_1^n)/\overline{h}(s_1^n)\right],
\end{aligned}
\end{equation}
where $\overline{f}(y_1^n)=\frac{1}{n}\sum_{j=1}^{n} \widehat{g}^{\beta-1}(y_j) f(y_j)$ and the last equality follows from the fact that $T(r_1^n)=T(s_1^n)$.
From the above expression, we note that $\mathcal{L}_{\mathcal{RE}}^{(\alpha,\beta)}(r_1 ^n;\lambda) - \mathcal{L}_{\mathcal{RE}}^{(\alpha,\beta)}(s_1^n;\lambda)$ is independent of $\lambda$. Thus, $T$ is a sufficient statistic for the $k$-parameter $\mathcal{M}^{(\alpha,\beta)}$-family with respect to the likelihood function $\mathcal{L}_{\mathcal{RE}}^{(\alpha,\beta)}$. $\hfill \Box$

From the above theorem, one may note that $\frac{\overline{f}}{\overline{h}}  = \left[\frac{\overline{f_1}}{\overline{h}},\dots ,\frac{\overline{f_d}}{\overline{h}}\right]^{T}$ is a sufficient statistic for the $k$-parameter $\mathcal{M}^{(\alpha,\beta)}$-family with respect to  $\mathcal{L}_{\mathcal{RE}}^{(\alpha,\beta)}$. Here we can observe that the sufficient statistic $\frac{\overline{f}}{\overline{h}}$ contains the elements $\frac{\overline{f_i}}{\overline{h}}$, for $1\leq i \leq d$, where $\overline{f_i}$ (resp. $\overline{h}$) is not just average of $f_i(y_j)$ (resp. $h(y_j)$), $j=1,\dots,n$ (as it was the case for the sufficient statistic corresponding $\mathcal{L^{(\alpha)}_J}$ derived in \cite{Generalized_sufficiency}).
In order to understand the significance of the presence of the tuning parameter $\beta$, let us consider $h(y) = 1$ and $f_i(y_j)=y_j$, for some $i$. Then, $\frac{\overline{f_i}}{\overline{h}}=\frac{\sum_{j=1}^n y_j\widehat{g}^{\beta-1}(y_j)}{\sum_{j=1}^n \widehat{g}^{\beta-1}(y_j)}$ is the weighted average of $y_j$ with coefficients $\frac{\widehat{g}^{\beta-1}(y_j)}{\sum_{j=1}^n\widehat{g}^{\beta-1}(y_j)}$. If we consider $\beta=1$, then $\frac{\overline{f_i}}{\overline{h}}=\frac{1}{n}\sum_{j=1}^n y_j$ which is the case for the estimation associated with $\mathcal{L^{(\alpha)}_J}$.  The presence of the tuning parameter $\beta$ in the coefficients can be helpful for robust estimation as it can control the weight of $y_j$ in order to reduce the effect of outliers. However, for $\beta=1$, equal weight is provided to each observation. We further explore the robustness of MLNREE (which is a function of sufficient statistic) in subsequent sections. Below we show that $\frac{\overline{f}}{\overline{h}} $ is indeed a minimal sufficient statistic for the $k$-parameter regular $\mathcal{M}^{(\alpha,\beta)}$-family with respect to $\mathcal{L}_{\mathcal{RE}}^{(\alpha,\beta)}$.
\begin{theorem}\label{T3} The statistic $ \frac{\overline{f}}{\overline{h}} = \left[\frac{\overline{f_1}}{\overline{h}}, \dots ,\frac{\overline{f_d}}{\overline{h}}\right]^T $, as defined in Theorem~\ref{T2}, is a minimal sufficient statistic for the $k$-parameter regular $\mathcal{M}^{(\alpha,\beta)}$-family with respect to $\mathcal{L}_{\mathcal{RE}}^{(\alpha,\beta)}$.
\end{theorem}
\noindent\textbf{Proof :} Let $r_1 ^n $ and $s_1^n$ be two samples such that $\left[\mathcal{L}_{\mathcal{RE}}^{(\alpha,\beta)}(r_1 ^n;\lambda) - \mathcal{L}_{\mathcal{RE}}^{(\alpha,\beta)}(s_1^n;\lambda)\right]$ is independent of $\lambda$. Now, from \eqref{E36}, we have that
\begin{equation*}
\mathcal{L}_{\mathcal{RE}}^{(\alpha,\beta)}(r_1 ^n;\lambda) - \mathcal{L}_{\mathcal{RE}}^{(\alpha,\beta)}(s_1^n;\lambda) =  \frac{1}{\alpha-\beta}\left[ \log \frac{\overline{h}(r_1^n)}{\overline{h}(s_1^n)}
+\log \frac{1+w(\lambda)^T\overline{f}(r_1^n)/\overline{h}(r_1^n)}{ 1+w(\lambda)^T\overline{f}(s_1^n)/\overline{h}(s_1^n)} \right]
\end{equation*}
is independent of $\lambda$. Let 
\begin{equation*}
  M=  \log \frac{1+w(\lambda)^T\overline{f}(r_1^n)/\overline{h}(r_1^n)}{ 1+w(\lambda)^T\overline{f}(s_1^n)/\overline{h}(s_1^n)}.
\end{equation*}
Then $M$ is independent of $\lambda$. Consequently, we have that

\begin{equation*}
    (1-\exp(M))+{w(\lambda)}^T\left [\overline{f}(r_1^n)/\overline{h}(r_1^n) - \exp(M)\left \{\overline{f}(s_1^n)/\overline{h}(s_1^n)\right \}\right] = 0.
\end{equation*}
Since the given family is regular, we have that $1, w_1, \dots , w_d$ are linearly independent. Consequently,  we get that $\exp(M) = 1$ and $\overline{f}(r_1^n)/\overline{h}(r_1^n) - \exp(M)\left\{\overline{f}(s_1^n)/\overline{h}(s_1^n)\right\} = 0.$ Further, these imply that $\overline{f}(r_1^n)/\overline{h}(r_1^n) = \left\{\overline{f}(s_1^n)/\overline{h}(s_1^n)\right\}$. Hence $\frac{\overline{f}}{\overline{h}}$ is a minimal sufficient statistic.$\hfill\Box$

\vspace{0.2cm}
Below we give an example wherein we obtain the minimal sufficient statistic for the family of Student's $t$-distributions with respect to $\mathcal{L}_{\mathcal{RE}}^{(\alpha,\beta)}$.
\begin{example}\label{Ex1}
Consider the Student's $t$-distribution as in Definition \ref{s_d}; its pdf is given by 
\begin{equation}\label{p_s_d}
  f_\lambda(y) = \frac{\Gamma(\frac{\nu+1}{2})}{\Gamma(\frac{\nu}{2})\sqrt{\pi \nu \sigma^2}} \left[1+\frac{1}{\nu} \left(\frac{y-\mu}{\sigma}\right)^2\right]^{-\frac{\nu+1}{2}},
  \end{equation}
where $\nu$ ($ > 1$) is the degrees of freedom  and  $\lambda = (\mu, \sigma^2) \in \R \times \R^+ $ are the parameters. Note that $f_{\lambda}$ can equivalently be written as
\begin{equation*}
  f_\lambda(y) =  N_\nu(\lambda)\left[1+ y^2 \frac{1}{\nu \sigma^2+\mu^2} - y \frac{2 \mu \sigma}{\nu \sigma^2+ \mu^2} \right]^{-\frac{\nu+1}{2}},
 \end{equation*}
where $ N_\nu(\lambda)=\frac{\Gamma(\frac{\nu+1}{2})}{\Gamma(\frac{\nu}{2})\sqrt{\pi \nu \sigma^2}} \left[1+\frac{ \mu^2}{\nu \sigma^2}\right]^{-\frac{2}{\nu+1}}$. Let $\alpha = \beta -\frac{2}{\nu+1}$ and $\frac{1}{\nu} = l_{\alpha, \beta} = \frac{\beta-\alpha}{2+\alpha-\beta}$. Then, the above pdf can equivalently be written as 
\begin{equation}\label{E37}
   f_\lambda(y) =  Z_{\alpha, \beta} (\lambda)\left[1+ y^2 \frac{l_{\alpha, \beta}}{ \sigma^2+ l_{\alpha, \beta}\mu^2} - y \frac{2 \mu \sigma l_{\alpha, \beta}}{ \sigma^2+ l_{\alpha, \beta}\mu^2} \right]^{\frac{1}{\alpha-\beta}},   
\end{equation}
where $Z_{\alpha, \beta} (\lambda)=N_{1/l_{\alpha, \beta}}(\lambda)$. Note that this forms a regular $\mathcal{M}^{(\alpha, \beta)}$-family with $h(y)=1$, $f(y) = [y^2, y]^T$ and $w(\lambda) = \left[\frac{l_{\alpha, \beta}}{ \sigma^2+ l_{\alpha, \beta}\mu^2}, -\frac{2 \mu \sigma l_{\alpha, \beta}}{ \sigma^2+ l_{\alpha, \beta}\mu^2}\right]^T$. By Theorem \ref{T2}, the sufficient statistic for this family of probability distribution is given by
$ \left[\frac{\sum_{j=1}^{n} y_j^{2}\widehat{g}^{\beta-1}(y_j)} {\sum_{j=1}^{n} \widehat{g}^{\beta-1}(y_j)}, \frac{\sum_{j=1}^{n} y_j \widehat{g}^{\beta-1}(y_j)} {\sum_{j=1}^{n} \widehat{g}^{\beta-1}(y_j)} \right]^T$. Again, by Theorem \ref{T3}, this is also a minimal sufficient statistic for the given family of Student's $t$-distributions.
\end{example}
\section{Generalized Rao-Blackwell theorem}\label{S5}
In this section, we derive the generalized Rao-Blackwell theorem for the minimum LNRE estimation which become useful to obtain the best robust estimator using a sufficient statistic. We begin this section by stating the classical 
Rao-Blackwell theorem following \cite{Blackwell,C.R. Rao}. 
\begin{theorem}
 Let $\mathcal{P}= \{f_\lambda: \lambda\in\Lambda\subset \R^k\}$ be a family of probability distributions, and let $\widehat{\lambda} $ be an unbiased estimator of $\lambda$. Further, let $T$ be a sufficient statistic for $\lambda$ with respect to the usual log-likelihood function $\mathcal{L}(y_1^n,\lambda)$ defined in \eqref{E8}. Then, $E_{\lambda}[\widehat{\lambda}|T]$ is independent of $\lambda$. Further, if $ \Phi(T) = \ E_{\lambda}[\widehat{\lambda}|T]$ then 
\begin{enumerate}
    \item [$(i)$]$E_{\lambda}[\Phi(T)] = E_{\lambda} [\widehat{\lambda}] = \lambda$;
    \item [$(ii)$]$Var_{\lambda} [\Phi(T)] \leq Var_{\lambda}[\widehat{\lambda}]$,
\end{enumerate}
for all $\lambda \in \Lambda$. Moreover, the equality in part ($ii$) holds if and only if $\widehat{\lambda}$ is a function of $T$. $\hfill\Box$
\end{theorem}
One may note that, for the exponential family defined in Definition~\ref{D2}, $\overline{f}= \frac{1}{n} \sum_{i=1}^{n} f(y_i)$ is a sufficient statistic for $\lambda$ with respect to the usual likelihood function. Let $\tau(\lambda)=E_{\lambda}[\overline{f}]=E_{\lambda}[f(Y)] $, which is the mean value parameter. Let $\widehat{\lambda}$ be any unbiased estimator of $\tau(\lambda)$. Then, by the Rao-Blackwell theorem, we have that
\begin{enumerate}
    \item [$(i)$]$\Phi(\overline{f}) := E_{\lambda} [\widehat{\lambda}|\overline{f}] = \tau(\lambda)$;

    \item [$(ii)$]$Var_{\lambda}[\Phi(\overline{f})] \leq Var_{\lambda}[\widehat{\lambda}]$.
\end{enumerate}
Under mild regularity conditions, it can be shown that there does not exist any unbiased estimator of $\tau(\lambda)$ that has lesser variance than that of $\overline {f}$. Moreover,

\begin{equation}\label{E38}
    Var_{\lambda}[\overline{f}] = \frac{[\tau'(\lambda)]^2}{n I(\lambda)},
\end{equation}
where $I(\lambda)= E_{\lambda} \left[ \left(\frac{\partial}{\partial \lambda} \log 
\hspace{0.05 in} f_{\lambda}(Y)\right)^2\right]$ is the Fisher information. Consequently, 
\begin{equation}\label{E39}
    Var_{\lambda}[\widehat{\lambda}] \geq \frac{[\tau'(\lambda)]^2}{n I(\lambda)},
\end{equation}
which is the famous  Cram\'{e}r-Rao inequality. In view of Fisher-Darmois-Koopman-Pitman theorem, one may note that the classical Rao-Blackwell theorem cannot be used outside the exponential family. Recently, Gayen and Kumar~\cite{Generalized_sufficiency} developed Rao-Blackwell type theorems for the minimum DPD and LDPD estimation. In what follows, we derive the generalized Rao-Blackwell theorem for the minimum LNRE estimation, starting with an important proposition.
\\\hspace*{0.2 in}Consider the likelihood function $\mathcal{L}_{\mathcal{RE}}^{(\alpha,\beta)}$ as defined in \eqref{Exx}. Then, the deformed probability distribution associated with $\mathcal{L}_{\mathcal{RE}}^{(\alpha,\beta)}$ is given by

\begin{equation}\label{E40}
\begin{aligned}
\widetilde f_{\lambda} (y_1^n)  := \frac{\exp\left[\mathcal{L}_{\mathcal{RE}}^{(\alpha,\beta)}(y_1^n;\lambda)\right]}{\bigints \exp \left[\mathcal{L}_{\mathcal{RE}}^{(\alpha,\beta)}(r_1 ^n;\lambda)\right] dr_1^n} 
= \frac{\left(\frac{1}{n}\sum_{j=1}^{n} \widehat{g}^{\beta-1}(y_j) f_\lambda ^{\alpha-\beta} (y_j)\right)^\frac{1}{\alpha-\beta}}{\bigints \left(\frac{1}{n}\sum_{j=1}^{n} \widehat{g}^{\beta-1}(r_j) f_\lambda ^{\alpha-\beta} (r_j)\right)^\frac{1}{\alpha-\beta}dr_1 ^n}. 
\end{aligned}
\end{equation}
By Proposition~\ref{GFD}, we can say that $T$ is a sufficient statistic for $\lambda$ with respect to  $\mathcal{L}_{\mathcal{RE}}^{(\alpha,\beta)}$ if and only if the associated conditional distribution $\widetilde{f}_{\lambda_{y_1^n|t}}(y_1^n|t)$ is independent of $\lambda$ for any value $t$ of $T$. Consequently, we have the following result.
\begin{proposition}\label{P2}
Let $\mathcal{P} = \{ f_\lambda : \lambda \in \Lambda \subset \R \}$ be a family of probability distributions with common support $S \subset \R$, and let $\widetilde{\mathcal{P}} = \{\widetilde{f}_\lambda : \lambda \in \Lambda \subset \R\}$ be the corresponding family of deformed probability distributions, where $\widetilde{f}_\lambda$ is as in \eqref{E40}. Let $\widehat{\lambda}$ be an unbiased estimator of $\widetilde{\tau}(\lambda) := \widetilde {E}_{\lambda}[\widehat{\lambda}]$. Further, let $T $ be a sufficient statistic for $\lambda$ with respect to $\mathcal{L}_{\mathcal{RE}}^{(\alpha,\beta)}$. Then, the conditional expectation $\widetilde{E}_\lambda[\widehat{\lambda}|T]$ is independent of $\lambda$. $\hfill\Box$
\end{proposition}

Since sufficient statistics play a crucial role in constructing estimators with lower variance, in the following theorem, we state a generalization of the Rao-Blackwell theorem for the minimum LNRE estimation.
\begin{theorem} \label{T4}Let $\mathcal{P} = \{ f_\lambda : \lambda \in \Lambda \subset \R \}$ be a family of probability distributions with common support $S \subset \R$, and let $\widetilde{\mathcal{P}} = \{\widetilde{f}_\lambda: \lambda\in\Lambda\subset\R\}$ be the corresponding family of deformed probability distributions, where $\widetilde{f}_\lambda$ is as defined in \eqref{E40}. Let $\widehat{\lambda}$ be an unbiased estimator for its expected value with respect to $\widetilde{f}_\lambda$, and let $\widetilde{\tau}(\lambda) = \widetilde{E}_\lambda[\widehat{\lambda}]$. Further, let $T$ be a sufficient statistic for $\lambda$ with respect to $\mathcal{L}_{\mathcal{RE}}^{(\alpha,\beta)}$, and let 
$\widetilde \varphi(T) = \widetilde{E}_\lambda[\widehat{\lambda}|T ]$. Then,
\begin{enumerate}
    \item [$(i)$] $\widetilde{E}_\lambda[\widetilde \varphi(T)] = \widetilde{\tau}(\lambda) $;
     \item [$(ii)$]$\widetilde{Var}_\lambda[\widetilde \varphi(T)] \leq \widetilde{Var}_\lambda[\widehat{\lambda}]$, 
\end{enumerate}
where 
the equality in part ($ii$) holds if and only if $\widehat{\lambda}$ is a function of $T$.
\end{theorem}
\noindent\textbf{Proof:} For a given partition set $C_{T(r_1^n)}$ of $S^n$, we have
\begin{eqnarray*}\label{E44}
        \widetilde \varphi(T(y_1^n)) & =& \int_{C_{T(y_1^n)}} \widehat{\lambda}(r_1^n) \widetilde{f}_{\lambda }(r_1^n|T(y_1^n)) d(r_1^n)\nonumber\\
        & =& \int_{C_{T(y_1^n)}} \widehat{\lambda}(r_1^n) \widetilde f_\lambda (r_1 ^n) d(r_1^n) \Bigg / \int_{C_{T(y_1^n)}} \widetilde f_\lambda (r_1 ^n)dr_1^n,
\end{eqnarray*}
and consequently, we get
\begin{equation*}
    \begin{aligned}
        \widetilde{E}_\lambda\left[\widetilde \varphi(T(Y_1^n)) \cdot \textbf{1}_{C_{T(r_1^n)}}\right] & = \int_{C_{T(r_1^n)}} \widetilde{f}_\lambda (y_1^n)\widetilde \varphi(T(y_1^n)) dy_1^n\\
        & = \frac{\bigintss_{C_{T(r_1^n)}} \widehat{\lambda}(y_1^n) \widetilde f_\lambda (y_1^n) d(y_1^n)} {\bigintss_{C_{T(r_1^n)}} \widetilde{f}_\lambda (y_1^n) dy_1^n} \left[\int_{C_{T(r_1^n)}} \widetilde{f}_\lambda (y_1^n) dy_1^n \right]\\
        & = \widetilde{E}_\lambda\left[\widehat{\lambda}(Y_1^n) \cdot \textbf{1}_{A_T(r_1^n)}\right].
    \end{aligned}
\end{equation*}
Note that the above equality holds for every partition set $C_{T(r_1^n)}$ of $S^n$. Thus, 
\begin{equation*}
  \widetilde{E}_\lambda[\widetilde \varphi(T(Y_1^n))] = \widetilde{E}_\lambda[\widehat{\lambda}(Y_1^n)] = \widetilde{\tau}(\lambda), 
\end{equation*}
and hence part (i) is proved. Further, note that
$$\widetilde{Var}_\lambda[\widetilde \varphi(T)] = \widetilde{E}_\lambda[\widetilde \varphi(T)]^2 - (\widetilde{E}_\lambda[\widetilde \varphi(T)])^2  \quad \text{and} \quad \widetilde{Var}_\lambda[\widehat{\lambda}] = \widetilde{E}_\lambda[\widehat{\lambda}]^2 -(\widetilde{E}_\lambda[\widehat{\lambda}])^2.$$
 Thus, to prove the result given in part (ii), it suffices to show that
$$\widetilde{E}_\lambda[\widetilde \varphi(T)]^2 \leq \widetilde{E}_\lambda[\widehat{\lambda}]^2,$$
or equivalently,
\begin{equation}\label{E*}
 \widetilde{E}_\lambda[\widetilde{E}_\lambda[\widehat{\lambda}|T ]]^2 \leq \widetilde{E}_\lambda[\widetilde{E}_\lambda[\widehat{\lambda}^2| T]].  
\end{equation}
Note that this holds if
$$\widetilde{E}^2_\lambda[\widehat{\lambda}|T ] \leq \widetilde{E}_\lambda[\widehat{\lambda}^2| T],$$
which follows from the Cauchy-Schwartz inequality given by
$$\widetilde{E}^2_\lambda[\widehat{\lambda}|T ] \leq \widetilde{E}_\lambda[\widehat{\lambda}^2| T] \widetilde{E}_\lambda[1| T].$$
Thus, $\widetilde{Var}_\lambda[\widetilde \varphi(T)] < \widetilde{Var}_\lambda[\widehat{\lambda}]$ is proved. Further, the equality given in part (ii) holds if and only if
$$\widetilde{E}_\lambda[\widetilde{E}_\lambda[\widehat{\lambda}^2| T]- (\widetilde{E}_\lambda[\widehat{\lambda}|T ])^2] = 0,$$
or equivalently,
$$\widetilde{E}_\lambda[\widetilde{Var}_\lambda[\widehat{\lambda}|T]] = 0.$$
Again, this is equivalent to the fact that 
$$\widetilde{E}_\lambda[\widehat{\lambda}^2| T] = \widetilde{E}^2_\lambda[\widehat{\lambda}|T ],$$
which holds if and only if $\widehat{\lambda}$ is a function of $T$. Hence, the result. $\hfill\Box$ 
\\\hspace*{0.2 in} From Theorem \ref{T4}, we can say that $\widetilde{\varphi}(T)$ is uniformly better estimator than $\widehat{\lambda}$ for $\widetilde{\tau}(\lambda)$. However, we can show that the best unbiased estimator of $\widetilde{\tau}(\lambda)$ with respect to $\widetilde{f}_\lambda$ is unique (independent of the choice of the sufficient statistic $T$) which is presented in the following. The proof follows along the same lines as in Lehmann and Casella~\cite{Lehmann & Casella_book} and is, therefore, omitted.
\begin{proposition}\label{P4} 
Let $\mathcal{P} = \{ f_\lambda : \lambda \in \Lambda \subset \R \}$ be a family of probability distributions with common support $S \subset \R$, and let $\widetilde{\mathcal{P}} = \{\widetilde{f}_\lambda: \lambda\in\Lambda\subset\R\}$ be the corresponding family of deformed probability distributions, where $\widetilde{f}_\lambda$ is as defined in \eqref{E40}.
Let $\widetilde{\varphi}_1(T)$ and $\widetilde{\varphi}_2(T)$ be any two unbiased estimators for $\widetilde{\tau}(\lambda)$ with respect to $\widetilde f_\lambda$. If the variances of $\widetilde{\varphi}_1(T)$ and $\widetilde{\varphi}_2(T)$ with respect to $\widetilde f_\lambda$ are equal and minimum, 
then $\widetilde{\varphi}_1(T)=\widetilde{\varphi}_2(T)$ almost surely.
\end{proposition}
\section{Generalized Cram\'{e}r-Rao bound}\label{S6}
In this section, we derive the generalized Cram\'{e}r-Rao bound for the variance of the estimators using the generalized Fisher information 
when the minimum LNRE estimation is under consideration. 
In particular, we derive  the variance of $\overline{f}/\overline{h}$,  a minimal sufficient statistic for regular $\mathcal{M}^{(\alpha,\beta)}$-family, and further investigate whether this variance is minimum or not.  

\begin{theorem}\label{T5}
Let $\mathcal{P} = \{ f_\lambda : \lambda \in \Lambda \subset \R \}$ be a family of probability distributions with common support $S \subset \R$, and let $\widetilde{\mathcal{P}} = \{\widetilde{f}_\lambda : \lambda \in \Lambda \subset \R\}$ be the corresponding family of deformed probability distributions, where $\widetilde{f}_\lambda$ is as in \eqref{E40}. Let $Y_1^{n}$ be an i.i.d. sample drawn from some $f_\lambda\in \mathcal{P}$, and let $y_1^n$ be its realized value. Suppose that the following conditions hold.
\begin{enumerate}
    \item [$(i)$] The derivative $\frac{\partial}{\partial\lambda} \log \hspace{0.05 in} \widetilde{f}_\lambda (y_1^n)$ exists and is finite;
    \item [$(ii)$] For any statistic $u$ with $\widetilde{E}_\lambda|u(Y_1^n)| < \infty$ for all $\lambda$, the operations of integration (summation) and differentiation with respect to $\lambda$ can be interchanged in $\widetilde{E}_\lambda|u(Y_1^n)| $;
\end{enumerate}
Let $\widehat{\lambda}$ be an unbiased estimator for its expected value with respect to $\widetilde{f}_\lambda$, and let $\widetilde{\tau}(\lambda) = \widetilde{E}_\lambda[\widehat{\lambda}(Y_1^n)]$. Then, 
\begin{equation}\label{E6.1}
 \widetilde{Var}_{\lambda}(\widehat{\lambda}(Y_1^n)) \geq \frac{\left[\frac{\partial}{\partial \lambda}{\widetilde{\tau}}(\lambda)\right]^2}{\widetilde{I}_n(\lambda)} ,  
\end{equation}
where $ \widetilde{I}_{n}(\lambda) = \widetilde{Var}_\lambda\left[\frac{\partial}{\partial\lambda} \log \hspace{0.05 in} \widetilde{f}_\lambda (Y_1^n)\right].$
\end{theorem}
\noindent \textbf{Proof:} Consider
\begin{eqnarray}\label{E6.2}
\widetilde{E}_{\lambda}\left[\frac{\partial}{\partial\lambda} \log \hspace{0.05 in} \widetilde{f}_\lambda (Y_1^n)\right] &=& \bigintsss \widetilde{f}_\lambda (y_1^n) \frac{\partial}{\partial\lambda} \log \hspace{0.05 in} \widetilde{f}_\lambda (y_1^n)dy_1^n \nonumber \\
&=& \frac{\partial}{\partial\lambda} \bigintsss \widetilde{f}_\lambda (y_1^n)dy_1^n = \frac{\partial}{\partial\lambda}(1) = 0.
\end{eqnarray}
Now, we have
$\widetilde{\tau}(\lambda) = \widetilde{E}_\lambda[\widehat{\lambda}(Y_1^n)].$
Differentiating this with respect to $\lambda$, we get
\begin{eqnarray*}
\frac{\partial}{\partial \lambda}{\widetilde{\tau}}(\lambda) = \int \widehat{\lambda}(y_1^n) \frac{\partial}{\partial \lambda} \widetilde{f}_\lambda(y_1^n) dy_1^n 
= \widetilde{Cov}_\lambda \left[\widehat{\lambda}(Y_1^n),  \frac{\partial}{\partial \lambda} \log \hspace{0.05 in} \widetilde{f}_\lambda(Y_1^n)\right],
\end{eqnarray*}
where the last equality follows from \eqref{E6.2}. Further, this, in view of the Cauchy-Schwarz inequality,  gives 
\begin{eqnarray*}
\left[\frac{\partial}{\partial \lambda}{\widetilde{\tau}}(\lambda)\right]^2 &=&  \widetilde{Cov}_\lambda ^2 \left[\widehat{\lambda}(Y_1^n),  \frac{\partial}{\partial \lambda} \log \hspace{0.05 in}\widetilde{f}_\lambda(Y_1^n)\right] \\
&\leq& \widetilde{Var}_{\lambda}(\widehat{\lambda}(Y_1^n)) \cdot \widetilde{Var}_\lambda\left[\frac{\partial}{\partial\lambda} \log \hspace{0.05 in} \widetilde{f}_\lambda (Y_1^n)\right]
\end{eqnarray*}
or equivalently, $$\widetilde{Var}_{\lambda}(\widehat{\lambda}(Y_1^n)) \geq \frac{\left[\frac{\partial}{\partial \lambda}{\widetilde{\tau}}(\lambda)\right]^2}{\widetilde{I}_n(\lambda)}.$$ Hence the result. 

\begin{remark}
The following observations can be made under the set-up and notation of Theorem \ref{T5}.
\begin{enumerate}
\item [$(a)$] For the ML estimation, the usual Fisher information is given by
$$ I(\lambda) = Var_\lambda \left[\frac{\partial}{\partial\lambda} \log \hspace{0.05 in} f_\lambda (Y_1^n)\right] = Var_\lambda\left[\frac{\partial}{\partial\lambda} \mathcal{L}(y_1^n, \lambda) \right],$$ which can be obtained from the generalized Fisher information
$$\widetilde{I}_{n}(\lambda) = \widetilde{Var}_\lambda\left[\frac{\partial}{\partial\lambda} \log \hspace{0.05 in} \widetilde{f}_\lambda (Y_1^n)\right] = \widetilde{Var}_\lambda\left[\frac{\partial}{\partial\lambda} \mathcal{L^{(\alpha, \beta)}_{RE}}(y_1^n, \lambda) \right]$$ by considering the usual log-likelihood function instead of $\mathcal{L^{(\alpha, \beta)}_{RE}}$;  
\item [$(b)$] For $\widetilde{\tau}(\lambda)=\lambda$,
$$\widetilde{Var}_{\lambda}(\widehat{\lambda}(Y_1^n)) \geq \frac{1}{\widetilde{I}_n(\lambda)}; $$
\item [$(c)$] From the inequality given in \eqref{E6.1}, we can say that, as the Fisher information increases, the lower bound for the variance of the estimator decreases and the best estimator will have smallest variance;
\item [$(d)$] If $\widetilde{I}_n(\lambda)$ (and $\frac{\partial}{\partial \lambda}{\widetilde{\tau}}(\lambda)$ is finite) or the variance of $\widehat{\lambda}$ is infinite, then the inequality given in \eqref{E6.1} trivially holds. $\hfill\Box$
\end{enumerate}
\end{remark}
\hspace*{0.02 in} In Theorem \ref{T3}, we showed that $\frac{\overline{f}}{\overline{h}}=\left[\frac{\overline{f_1}}{\overline{h}},\dots ,\frac{\overline{f_d}}{\overline{h}}\right]^T$ is a minimal sufficient statistic for the regular $\mathcal{M}^{(\alpha,\beta)}$-family.  Thus, one may be interested to know whether this attains the Cram\'{e}r-Rao lower bound given in \eqref{E6.1}. To answer this question, we first obtain the variance of ${\overline{f}}/{\overline{h}}$ in the following theorem.
\begin{theorem}\label{T6}
Consider the regular $\mathcal{M}^{(\alpha,\beta)}$-family as in \eqref{E17}. Let $\widetilde \tau(\lambda) = \widetilde{E}_\lambda[\overline{f} / \overline{h}]$. Then the variance of the estimator $\overline{f} / \overline{h}$ is given by
\begin{equation*}
\widetilde{Var}_\lambda[\overline{f}/ \overline{h} ] = \frac{\left[\frac{\partial}{\partial \lambda} \widetilde{\tau}(\lambda)\right]^2}{\widetilde{I}_{\alpha,n} ^\beta (\lambda)},
\end{equation*}
where 
\begin{equation}\label{E54}
\widetilde{I}_{\alpha,n} ^\beta (\lambda) = \frac{\left\{\widetilde{Cov}_\lambda[\widetilde{s}(Y_1^n,\lambda), {\widetilde{f}_\lambda (Y_1^n)}^{\alpha-\beta} \widetilde{s}(Y_1^n,\lambda) / \overline{h}] \right\}^2}{\widetilde{Var}_\lambda[{\widetilde{f}_\lambda (Y_1^n)}^{\alpha-\beta} \widetilde{s}(Y_1^n,\lambda) / \overline{h}]} \quad \text{ and } \quad \widetilde{s}(y_1^n,\lambda) = \frac{\partial}{\partial\lambda} \log \hspace{0.05 in} \widetilde{f}_\lambda (y_1^n). 
\end{equation}
\end{theorem} 
\noindent\textbf{Proof :} 
Let $f_\lambda \in \mathcal{M}^{(\alpha,\beta)}$. Then, 
\begin{equation*}
    f_\lambda (y) = N(\lambda){[h(y)+w(\lambda) f(y)]}^\frac{1}{\alpha-\beta},\quad y\in S.
\end{equation*}
Consequently, in view of \eqref{E40}, the corresponding deformed pdf  associated with $\mathcal{L}_{\mathcal{RE}}^{(\alpha,\beta)}$ is given by
\begin{eqnarray}
\widetilde f_{\lambda} (y_1^n) 
 = [G(\lambda)\overline{h}+ \widetilde w(\lambda)\overline{f}]^{\frac{1}{\alpha - \beta}},\label{E47}
\end{eqnarray}
where $G(\lambda)=\left(\frac{N(\lambda)}{Z(\lambda)}\right)^{\alpha-\beta}$, $Z(\lambda)=\bigints \left(\frac{1}{n}\sum_{j=1}^{n} \widehat{g}^{\beta-1}(r_j) N(\lambda)^{\alpha - \beta}[h(r_j)+w(\lambda)f(r_j)]\right)^\frac{1}{\alpha-\beta}dr_1 ^n$,  $\widetilde w(\lambda)=w(\lambda)G(\lambda), \overline{f} = \frac{1}{n}\sum_{j=1}^{n} \widehat{g}^{\beta-1}(y_j) f(y_j)$ and $\overline{h} = \frac{1}{n}\sum_{j=1}^{n} \widehat{g}^{\beta-1}(y_j) h(y_j)$. Differentiating \eqref{E47} with respect to $\lambda$,  we get 
\begin{equation*}\label{E56}
    \begin{aligned}
 \frac{\partial}{\partial \lambda}\widetilde f_{\lambda} (y_1^n)
       = \frac{1}{\alpha - \beta} \widetilde f_{\lambda} (y_1^n)^ {1-\alpha+\beta}[G'(\lambda)\overline{h}+ \widetilde w'(\lambda)\overline{f}].
    \end{aligned}
\end{equation*}
Now,
\allowdisplaybreaks{
    \begin{eqnarray}
&&\widetilde{Cov}_\lambda[\widetilde{s}(Y_1^n,\lambda), {\widetilde{f}_\lambda (Y_1^n)}^{\alpha-\beta} \widetilde{s}(Y_1^n,\lambda) / \overline{h}] \nonumber \\ 
&& = \widetilde E_\lambda[\widetilde{s}^{2}(Y_1^n,\lambda) {\widetilde{f}_\lambda (Y_1^n)}^{\alpha-\beta} / \overline{h}] \nonumber\\ 
&& = \int \frac{\widetilde{s}(y_1^n,\lambda)}{\overline{h}} \frac{1}{\widetilde{f}_\lambda (y_1^n)}\left[\frac{\partial}{\partial\lambda} \widetilde{f}_\lambda (y_1^n) \right] {\widetilde{f}_\lambda (y_1^n)}^{\alpha-\beta+1}   dy_1^n \nonumber\\
&& = \frac{1}{\alpha - \beta} \int \frac{\widetilde{s}(y_1^n,\lambda)}{\overline{h}} \widetilde{f}_\lambda (y_1^n)[G'(\lambda)\overline{h}+ \widetilde w'(\lambda)\overline{f}]   dy_1^n \nonumber\\
&& = \frac{\widetilde w'(\lambda)}{\alpha - \beta} \int \left[\frac{\partial}{\partial\lambda} \widetilde{f}_\lambda (y_1^n) \right]\frac{   \overline{f}}{\overline{h}} dy_1^n \nonumber\\
&& = \left (\frac{\widetilde w'(\lambda)}{\alpha - \beta} \right) \frac{\partial}{\partial\lambda}\widetilde \tau(\lambda) .\label{E57}
\end{eqnarray}}
and
\begin{align}
\widetilde{Var}_\lambda[{\widetilde{f}_\lambda (Y_1^n)}^{\alpha-\beta} \widetilde{s}  (Y_1^n,\lambda) /  \overline{h}] 
& = \widetilde E_\lambda[\{{\widetilde{f}_\lambda (Y_1^n)}^{\alpha-\beta}\widetilde{s}(Y_1^n,\lambda)  / \overline{h}\}^2] - (\widetilde E_\lambda[{\widetilde{f}_\lambda (Y_1^n)}^{\alpha-\beta} \widetilde{s}(Y_1^n,\lambda) / \overline{h}])^2 \nonumber\\
& = \left(\frac{1}{\alpha - \beta}\right)^2 \widetilde{Var}_\lambda[G'(\lambda)+ \widetilde w'(\lambda)(\overline{f}/{\overline{h})}] \nonumber\\
& = \left(\frac{1}{\alpha - \beta}\right)^2 [\widetilde w'(\lambda)]^{2}\widetilde{Var}_\lambda[\overline{f}/{\overline{h}}]. \label{E58}
\end{align}
By using \eqref{E57} and \eqref{E58}, we get
\begin{equation*}
        \frac{1}{\widetilde{I}_{\alpha,n} ^\beta (\lambda)} = \frac{\left(\frac{1}{\alpha - \beta}\right)^2 [\widetilde w'(\lambda)]^{2}\widetilde{Var}_\lambda[\overline{f}/{\overline{h}}]}{\left[\left(\frac{\widetilde w'(\lambda)}{\alpha - \beta} \right)\frac{\partial}{\partial\lambda}\widetilde \tau(\lambda)\right]^2},
\end{equation*}
or equivalently,
\begin{equation}\label{E59}
    \widetilde{Var}_\lambda[\overline{f}/{\overline{h}}] =  \frac{{\left[ \frac{\partial}{\partial\lambda}\widetilde \tau(\lambda)\right]^2}}{\widetilde{I}_{\alpha,n} ^\beta (\lambda)}.
\end{equation}
Hence the result.
\begin{remark}\label{R4}
The following observations can be made under the set-up and notation of Theorem \ref{T6}.
\begin{enumerate}
    \item [$(a)$] For $h= 1$, $\alpha \rightarrow 1 $ and $\beta \rightarrow 1$, $\widetilde{I}_{\alpha,n} ^\beta (\lambda)$ reduces to   $ \widetilde{I}_{n}(\lambda)$;
     \item [$(b)$] By Theorems~\ref{T5} and \ref{T6}, we have 
$$\widetilde{Var}_\lambda[\overline{f}/{\overline{h}}] \geq \frac{\left[\frac{\partial}{\partial \lambda}{\widetilde{\tau}}(\lambda)\right]^2}{\widetilde{I}_n(\lambda)}. $$
Moreover, the equality in the above expression holds when $h= 1$, $\alpha \rightarrow 1 $ and $\beta \rightarrow 1$. $\hfill\Box$
\end{enumerate} 
\end{remark}  
\hspace*{0.02 in}In the previous remark, we see that $\overline{f}/\overline{h}$ does not attain the Cram\'{e}r-Rao  lower bound, apart from a particular case. Thus, one immediate question arises whether
$\overline{f}/\overline{h}$ is the best unbiased estimator for its expected value for the regular $\mathcal{M^{(\alpha, \beta)}}$-family. Interestingly, this is not the case. Below we give a counterexample wherein we show that there exists an unbiased estimator whose variance is less than that of $\overline{f}/\overline{h}$.
\begin{example}\label{Ex3}
 Consider the family of Bernoulli distributions with the pdf given by
 \begin{equation*}
    f_\lambda(y) = (1-\lambda)\left[1+ y\frac{2\lambda-1}{1-\lambda}\right], \quad y=0,1,
 \end{equation*}
 where $\lambda \in (1/2,1)$. Note that this belongs to the regular $\mathcal{M^{(\alpha, \beta)}}$-family when $\alpha=2$, $\beta=1$, $N(\lambda) = 1-\lambda$, $w(\lambda) = \frac{2\lambda-1}{1-\lambda}$, $f(y)= y$ and $h(y) = 1$. Let $Y_1 ,Y_2$ be a random sample drawn according to $f_\lambda$. Then, the deformed pdf associated with $\mathcal{L}_{\mathcal{RE}}^{(\alpha,\beta)}$ is given by
 \begin{equation*}
     \widetilde{f}_{\lambda}(y_1, y_2) = \frac{1-\lambda}{2}\left[1+\overline{y}\frac{2\lambda-1}{1-\lambda}\right],
 \end{equation*}
where $\overline{y}= \frac{\overline f}{\overline h}$, $\overline f = \frac{y_1 + y_2}{2} $ and $\overline h = 1$. Let $\widetilde{E}_\lambda[\overline{Y}] = \widetilde \tau(\lambda) $. To prove that there exists other unbiased estimator of $\widetilde \tau(\lambda)$ whose variance is less than the variance of $\overline{Y}$, consider a function $h$ of random vector $(Y_1, Y_2)$ as
\begin{equation*}
    h(Y_1, Y_2) = \begin{cases}
 1 , &\text{if $(Y_1, Y_2)= (0,0)$ or $(1,1)$}\\
-1, &\text{if $(Y_1, Y_2)= (0,1)$ or $(1,0)$}.
    \end{cases}
\end{equation*}
Then, $\widetilde{E}_\lambda[h(Y_1, Y_2)] = 0$. Now, consider another function $k$ of random vector $(Y_1, Y_2)$ as $k(Y_1, Y_2) = \overline{Y} - \frac{1}{8} h(Y_1, Y_2)$. Then, $k(Y_1, Y_2)$ is also an unbiased estimator of $\widetilde \tau(\lambda)$. Next, we proceed to show that $\widetilde{Var}_{{3}/{4}} (\overline{Y}) > \widetilde{Var}_{{3}/{4}} (k(Y_1, Y_2)).$ 
Now, we have
\begin{equation*}
\widetilde{E}_{3/4}[\overline{Y}]^2 = \sum_{(y_1, y_2)} { \overline{y} }^2   \widetilde{f}_{3/4}(y_1, y_2) = 0.5
\end{equation*}
and
\begin{equation*}
 \widetilde{E}_{3/4}[k(Y_1, Y_2)]^2 = \sum_{(y_1, y_2)} [\overline{y}-{1}/{8}h(y_1, y_2)]^2   \widetilde{f}_{3/4}(y_1, y_2) = 0.4843.   
\end{equation*}
Consequently, we get
\begin{equation*}
 \widetilde{E}_{3/4}[\overline{Y}]^2 > \widetilde{E}_{3/4}[k(Y_1, Y_2)]^2 ,  
\end{equation*}
or equivalently,
\begin{equation*}
 \widetilde{Var}_{{3}/{4}} (\overline{Y}) > \widetilde{Var}_{{3}/{4}} (k(Y_1, Y_2)).  
\end{equation*}
Thus, the variance of $\overline{Y}$ is strictly greater than that of $k(Y_1, Y_2)$. Hence $\overline{Y}$ is not the best estimator for its expected value with respect to the generalized likelihood function associated with the LNRE.
\end{example}
\begin{remark}\label{R5}
 Note that $\mathcal{M}^{(\alpha, \beta)}$-family coincides with the $\mathcal{M}^{(\alpha)}$-family when $\beta=1$.
    Thus, in light of the previous example, we may conclude that the minimal sufficient statistic, $ T = \frac{\overline{f}}{\overline{h}} = \left[\frac{\overline{f_1}}{\overline{h}}, \dots ,\frac{\overline{f_d}}{\overline{h}}\right]^T $, of the regular $\mathcal{M^{(\alpha)}}$-family with $w(\lambda)>0$ is not the best estimator for its expected value. However, this contradicts with the statement given in Theorem $20$ of \cite{Generalized_sufficiency}.
   This contradiction arises due to a loophole in proving the covariance to be greater than or equal to zero. In the proof of the said theorem, they obtained that, for $k(\theta)\neq 0$,
\begin{equation}\label{Corr1}
(\alpha-1) \frac{\partial}{\partial \theta}\left[\frac{k(\theta)}{N(\theta)}\right]\Big/\left[\frac{k(\theta)}{N(\theta)}\right] = \frac{w'(\theta)}{w(\theta)},
\end{equation}
and claimed that this expression \eqref{Corr1} can equivalently be written as
\begin{equation}\label{Corr2}
\frac{\partial}{\partial \theta} \log \left[\frac{k(\theta)}{N(\theta)}\right]^{\alpha-1} = \frac{\partial}{\partial \theta}\log{w(\theta)},  \end{equation}
which is not true in general. To obtain \eqref{Corr2} from  \eqref{Corr1}, we must have the condition that $k(\theta) > 0$ (since $N(\theta) > 0$), which is needed to be proved. This is because $\log(\cdot)$ function is defined only for positive arguments.
\end{remark}
\section{The MLNREE for the family of Student's distributions}\label{S7}
 In this section, we first obtain the MLNREE for the family of Student's $t$-distributions. Note that this estimator is  a function of sufficient statistic as mentioned in Example \ref{Ex1}. Then, we do a simulation study  in order to highlight the utility of sufficient statistic as well as robustness and flexibility of this estimator. Next, we perform a similar study for the family of Student's $r$-distributions. Moreover, we  discuss a real data where we compare the performance of the MLNREE with that of other related estimators. 
\subsection{Student's $t$-distribution with positive degrees of freedom }\label{S7.1}
Let $Y_1, \dots , Y_n$ be an i.i.d sample from a true distribution $g$, which is a Student's $t$-distribution with the pdf given by \eqref{p_s_d},
where $\nu \in (1, \infty)$ is assumed to be known, and $\lambda = (\mu, \sigma^2) \in \Lambda = \R\times\R^+ \subset \R^2 $ is the set of unknown parameters. Note that $f_{\lambda}(y) \in \mathcal{M}^{(\alpha, \beta)}$- family, for $\beta = \alpha + \frac{2}{\nu+1} $. Further, we have 
\begin{equation}\label{ar01}
    u(y, \lambda ) = \nabla \log \hspace{0.05 in} f_\lambda (y) = \left[\frac{1}{f_{\lambda}(y)} \frac{\partial f_{\lambda}(y)}{\partial \mu}, \frac{1}{f_{\lambda}(y)} \frac{\partial f_{\lambda}(y)}{\partial \sigma^2}\right]  = \left[\frac{(\nu+1)(y-\mu)}{\nu \sigma^2 \left[1+\frac{(y-\mu)^2}{\nu \sigma^2}\right]}, \frac{(y-\mu)^2 - \sigma^2}{2\sigma^4\left[1+\frac{(y-\mu)^2}{\nu \sigma^2}\right]}\right]. 
\end{equation}
In order to obtain the estimators of $\lambda=(\mu, \sigma^2 )$, we maximize  $\mathcal{L}_{\mathcal{RE}}^{(\alpha,\beta)}$ given in \eqref{E16}, for $\mu$ and $\sigma^2$. Now, the second term of $\mathcal{L}_{\mathcal{RE}}^{(\alpha,\beta)}$ is given by
\begin{equation*}
\begin{aligned}
\frac{1}{\alpha}\log \int f_\lambda^\alpha(y) dy & = \frac{1}{\alpha}\log \left(\frac{\Gamma(\frac{\nu+1}{2})}{\Gamma(\frac{\nu}{2})\sqrt{\pi \nu \sigma^2}} \right)^\alpha \bigintsss_{-\infty}^{+\infty} {1}/{\left[1+\frac{(y-\mu)^2}{\nu \sigma^2}\right]^{\alpha(\frac{\nu+1}{2})}}  dy\\
& = \frac{1}{\alpha}\log \left(\frac{\Gamma(\frac{\nu+1}{2})}{\Gamma(\frac{\nu}{2})\sqrt{\pi \nu \sigma^2}} \right)^\alpha \bigintsss_{-\infty}^{+\infty} {\sigma}/{\left[1+\frac{t^2}{\nu}\right]^{\alpha(\frac{\nu+1}{2})}}  dt,\\
\end{aligned}
\end{equation*}
where the second equality is obtained by applying the transformation $\frac{y-\mu}{\sigma} = t$. This final result  becomes independent of $\mu$. Thus, to obtain the estimator of $\mu$, it suffices to consider only the first term of $\mathcal{L}_{\mathcal{RE}}^{(\alpha,\beta)}$. Differentiating $\mathcal{L}_{\mathcal{RE}}^{(\alpha,\beta)}$ with respect to $\mu$ and then equating  it to zero, we get 
\begin{equation*}
\frac{\frac{1}{n}\sum_{j=1} ^ {n} \widehat{g}^{\beta-1}(y_j) f_\lambda^{\alpha-\beta}(y_j) \frac{1}{f_{\lambda}(y_j)} \frac{\partial f_{\lambda}(y_j)}{\partial \mu}} {\frac{1}{n}\sum_{j=1} ^ {n} \widehat{g}^{\beta-1}(y_j) f_\lambda^{\alpha-\beta}(y_j)} = 0.
\end{equation*}
By using \eqref{p_s_d} and \eqref{ar01} in the above equation, we get
\begin{equation}\label{mu_RE}
\widehat{\mu}_{\mathcal{RE}} = \frac{\sum_{j=1}^{n} y_j \widehat{g}^{\beta-1}(y_j)} {\sum_{j=1}^{n} \widehat{g}^{\beta-1}(y_j)} = T_1 \quad \text{(say)}.
\end{equation}
Next, we proceed to obtain the estimator of $\sigma^2$. Differentiating $\mathcal{L}_{\mathcal{RE}}^{(\alpha,\beta)}$ with respect to $\sigma^2$ and equating it to zero, we get 
\begin{equation*}
\frac{\frac{1}{n} \sum_{j=1} ^ {n} \widehat{g}^{\beta-1}(y_j) f_{\lambda}^{\alpha-\beta}(y_j) \frac{1}{f_{\lambda}(y_j)} \frac{\partial f_{\lambda}(y_j)}{\partial \sigma^2}}{\frac{1}{n} \sum_{j=1} ^ {n} \widehat{g}^{\beta-1}(y_j) f_{\lambda}^{\alpha-\beta}(y_j)} = \frac{\int f_{\lambda}^\alpha(y)\frac{1}{f_{\lambda}(y)} \frac{\partial f_{\lambda}(y)}{\partial \sigma^2}  dy}{\int f_{\lambda}^\alpha(y)  dy}. 
\end{equation*}
By using \eqref{p_s_d} and \eqref{ar01} in the above equation, we get
\begin{equation*}
 \frac{ \sum_{j=1} ^ {n} \widehat{g}^{\beta-1}(y_j) [{(y_{j}-\mu)}^2 - \sigma^2]}{ \sum_{j=1} ^ {n} \widehat{g}^{\beta-1}(y_j) \left[1+\frac{(y_j-\mu)^2}{\nu \sigma^2}\right]} = \frac{\bigints_{-\infty}^{+\infty} [(y-\mu)^2 - \sigma^2]/{\left[1+\frac{(y-\mu)^2}{\nu \sigma^2}\right]^{1+\alpha(\frac{\nu+1}{2})}}  dy}{\bigints_{-\infty}^{+\infty} {1}/{\left[1+\frac{(y-\mu)^2}{\nu \sigma^2}\right]^{\alpha(\frac{\nu+1}{2})}}  dy}.    
\end{equation*}
Since $\alpha = \beta-\frac{2}{\nu+1}$, the above equality can equivalently be written as 
\begin{equation}\label{ar03}
 \frac{ \sum_{j=1} ^ {n} \widehat{g}^{\beta-1}(y_j) [{(y_{j}-\mu)}^2 - \sigma^2]}{ \sum_{j=1} ^ {n} \widehat{g}^{\beta-1}(y_j) [\nu \sigma^2+(y_j-\mu)^2]} =C_1,       
\end{equation}
{where 
$$C_1=\frac{\bigintsss_{-\infty}^{+\infty} [t^2-1]/{[\nu+t^2]^{\frac{(\nu+1)\beta}{2}}}  dt}{\bigintsss_{-\infty}^{+\infty} {1}/{[\nu+t^2]^{\frac{(\nu+1)\beta}{2}-1}}  dt}= 1-\frac{\frac{\nu+1}{\nu}\left[\left(\frac{\nu+1}{2}\right)\beta-\frac{3}{2}\right]}{\left[\left(\frac{\nu+1}{2}\right)\beta-1\right]},$$
provided $\beta > \frac{3}{\nu+1}$}. Now, solving Equation \eqref{ar03} for $\sigma^2$, we get the MLNREE of $\sigma^2$ as
\begin{eqnarray}
    \widehat{\sigma}_{\mathcal{RE}}^2 = \left[\frac{1-C_1}{1+\nu C_1}\right]\frac{\sum_{j=1} ^ {n} {(y_{j}-\widehat{\mu}_{\mathcal{RE}})}^2\widehat{g}^{\beta-1}(y_j)}{\sum_{j=1} ^ {n}\widehat{g}^{\beta-1}(y_j)}=\left[\frac{1-C_1}{1+\nu C_1}\right] (T_2-T_1^2),\label{E60}
\end{eqnarray}
where
$$T_2=\frac{\sum_{j=1} ^ {n} {y_{j}}^2\widehat{g}^{\beta-1}(y_j)}{\sum_{j=1} ^ {n}\widehat{g}^{\beta-1}(y_j)}.$$ 
Thus, the MLNREE of $\lambda=(\mu$, $\sigma^2)$ is given by $\widehat{\lambda}=(\widehat{\mu}_{\mathcal{RE}},\widehat{\sigma}_{\mathcal{RE}}^2 )$, for $\beta > \frac{3}{\nu+1}$, where $\widehat{\mu}_{\mathcal{RE}}$ and $\widehat{\sigma}_{\mathcal{RE}}^2$ are given in \eqref{mu_RE} and \eqref{E60}, respectively. Moreover, one may observe that both the estimators are functions of sufficient statistic for the family of Student's $t$-distributions with $\nu > 1$ as discussed in Example \ref{Ex1}.
\\\hspace*{0.2 in}As mentioned in Gayen and Kumar~\cite{Projection_theorems},  the MDPDE and MLDPDE for $\lambda=(\mu, \sigma^2)$  are the same, and are given by 
\begin{equation}
 \widehat{\mu}_{\mathcal{B}} = \widehat{\mu}_{\mathcal{L}} = \frac{1}{n} \sum_{j=1}^n y_j     
\end{equation}
and
\begin{equation}\label{E61}
    \widehat{\sigma}_{\mathcal{B}}^2 = \widehat{\sigma}_{\mathcal{L}}^2 = \left(\frac{\nu-2}{\nu}\right)\frac{1}{n} \sum_{j=1}^n {(y_j-\widehat{\mu}_{\mathcal{B}})}^2,
\end{equation}
which also follow from \eqref{mu_RE} and \eqref{E60}, respectively, for $\beta=1$. In case of minimum DPD and minimum LDPD estimations, one may note that the tuning parameter $\alpha$ and the degrees of freedom $\nu$ are related  as $\alpha = 1-\frac{2}{\nu+1}$. Thus, for a fixed $\nu$, the tuning parameter $\alpha$ is also fixed and consequently, we get unique MDPDE and MLDPDE. On the other hand, MLNREEs are not unique; rather we get an MLNREE for each choice of $\beta$ (another tuning parameter).  Moreover, by judiciously choosing $\beta$, we can obtain the best robust MLNREE which is justified below through a simulation study by considering different degrees of contamination in data.
\\\hspace*{0.2 in}We first  draw  a sample of size $n= 50$ from the mixture $(1-\eta)g+\eta N(0,16)$, where $g$ is the pdf of Student's $t$-distribution  with $\nu=3$, $\mu=0$ and $\sigma^2=1$, and $\eta$ represents the contamination level. In our study, we consider $\eta=0.0, 0.05, 0.10, 0.15, 0.20$ and $0.25$.
Then, based on the simulated sample, we compute the MLNREE, i.e., $\widehat{\lambda}=(\widehat{\mu}_{\mathcal{RE}}, \widehat{\sigma}_{\mathcal{RE}}^2)$ for $\lambda=(\mu, \sigma^2)$ for several choices of $(\alpha, \beta)$ such that $\alpha = \beta-\frac{1}{2} $ ($>0$) and $\beta>\frac{3}{4}$ hold, i.e., $\beta \in (3/4, \infty)$. Here, the non-parametric continuous estimate of density from the sample is taken to be the kernel density estimate with Gaussian kernel. We replicate the process $1000$ times and report the average value of the MLNREE along with their standard errors (SE) in Table \ref{Table_1}.

In Table \ref{Table_1}, the MLNREEs of $\sigma^2$ corresponding to different values of $\beta$ are marked in bold font, and the corresponding MLNREEs of $\mu$ are also given. Further, the MDPDEs and MLDPDEs are obtained corresponding to $\beta=1$. 
In Table \ref{yu}, we summarize the MDPDEs and MLDPDEs of $(\mu,\sigma^2)$ as well as the best MLNREE of $\sigma^2$ (along with corresponding estimate of $\mu$), for different contaminated data. From this table, we conclude that the MDPDE/ MLDPDE of $\sigma^2$ goes far away from the true value (i.e., $\sigma^2=1$) as the degree of contamination increases. Due to the flexibility of choosing the tuning parameter $\beta$, the best MLNREE of $\sigma^2$ is always close to the true parameter value irrespective of the increased degree of contamination. Moreover, one may note that the MLNREEs of $\mu$ remain near to the true value for all values of $\beta$. We can also see that, as the degree of contamination increases, the value of $\beta$ decreases in order to get better estimate.  Thus, we may conclude that, for estimating $\sigma^2$, the MLNREE  outperforms the MDPDE and MLDPDE for data contaminated by outliers at least for our empirical settings.
    \renewcommand{\arraystretch}{1.5}
 \begin{longtable}{|>{\centering\arraybackslash}p{0.8cm}|>{\centering\arraybackslash}p{2cm}|>{\centering\arraybackslash}p{2cm}|>{\centering\arraybackslash}p{2cm}|>{\centering\arraybackslash}p{2cm}|>{\centering\arraybackslash}p{2cm}|>{\centering\arraybackslash}p{2cm}|>{\centering\arraybackslash}p{2cm}|}
 \caption{Average MLNREEs corresponding to different choices of $\beta$, for  different contamination levels, along with their SEs in the parenthesis.}\label{Table_1} \\
    \hline
    $\beta$ & $0.85$ & $0.86$ & $0.88$ & $0.91$ & $0.96$ & $1$ & $1.02$ \\
    \hline
    
    
    
         \multicolumn{8}{|c|}{For samples without contamination. } \\
        \hline
        \makecell{$\widehat{\mu}_{\mathcal{RE}}$\\{\scriptsize$(SE)$}} & \makecell{$-0.015$\\{\scriptsize$(8.18\times10^{-18})$}}  & \makecell{$-0.015$\\{\scriptsize$(1.30\times10^{-18})$}}  & \makecell{$-0.015$\\{\scriptsize$(9.82\times10^{-18})$}}  & \makecell{$-0.014$\\{\scriptsize$(2.74\times10^{-18})$}}  & \makecell{$-0.013$\\{\scriptsize$(9.02\times10^{-19})$}}  & \makecell{$-0.012$\\{\scriptsize$(2.05\times10^{-18})$}}  & \makecell{$-0.012$\\{\scriptsize$(4.79\times10^{-19})$}} \\
        \hline
        \makecell{$\widehat{\sigma}^2_{\mathcal{RE}}$\\{\scriptsize$(SE)$}}  & \makecell{$0.487$\\{\scriptsize$(1.52\times10^{-17})$}} & \makecell{$0.527$\\{\scriptsize$(1.57\times10^{-17})$}} & \makecell{$0.604$\\{\scriptsize$(2.04\times10^{-17})$}} & \makecell{$0.711$\\{\scriptsize$(2.79\times10^{-17})$}} & \makecell{$0.864$\\{\scriptsize$(1.63\times10^{-17})$}} & \makecell{$0.968$\\{\scriptsize$(1.05\times10^{-17})$}} & \makecell{$\textbf{1.014}$\\{\scriptsize$(1.43\times10^{-16})$}} \\
        \hline
        \multicolumn{8}{|c|}{For samples with $5\%$ contamination $(0.95g+0.05N(0,16))$.} \\
        \hline
        \makecell{$\widehat{\mu}_{\mathcal{RE}}$\\{\scriptsize$(SE)$}} & \makecell{$-0.004$\\{\scriptsize$(5.88\times10^{-18})$}} & \makecell{$-0.004$\\{\scriptsize$(2.77\times10^{-19})$}} &  \makecell{$-0.004$\\{\scriptsize$(2.77\times10^{-19})$}} &  \makecell{$-0.003$\\{\scriptsize$(1.13\times10^{-17})$}}& \makecell{$-0.003$\\{\scriptsize$(1.27\times10^{-18})$}} & \makecell{$-0.003$\\{\scriptsize$(5.77\times10^{-18})$}} & \makecell{$-0.003$\\{\scriptsize$(8.27\times10^{-18})$}}\\
        \hline
         \makecell{$\widehat{\sigma}^2_{\mathcal{RE}}$\\{\scriptsize$(SE)$}} & \makecell{$0.574$\\{\scriptsize$(1.16\times10^{-17})$}} & \makecell{$0.622$\\{\scriptsize$(2.94\times10^{-17})$}} & \makecell{$0.713$\\{\scriptsize$(6.23\times10^{-17})$}} & \makecell{$0.838$\\{\scriptsize$(8.15\times10^{-17})$}} & \makecell{$\textbf{1.019}$\\{\scriptsize$(4.10\times10^{-17})$}} & \makecell{$1.141$\\{\scriptsize$(1.10\times10^{-16})$}} & \makecell{$1.195$\\{\scriptsize$(6.46\times10^{-17})$}} \\
        \hline
        \multicolumn{8}{|c|}{For samples with $10\%$ contamination $(0.90g+0.10N(0,16))$.} \\
        \hline
        \makecell{$\widehat{\mu}_{\mathcal{RE}}$\\{\scriptsize$(SE)$}} & \makecell{$-0.002$\\{\scriptsize$(2.22\times10^{-19})$}} & \makecell{$-0.002$\\{\scriptsize$(8.94\times10^{-18})$}} & \makecell{$-0.001$\\{\scriptsize$(2.38\times10^{-18})$}} & \makecell{$-0.001$\\{\scriptsize$(4.38\times10^{-18})$}} & \makecell{$-0.001$\\{\scriptsize$(9.44\times10^{-19})$}} & \makecell{$-0.001$\\{\scriptsize$(8.24\times10^{-18})$}} & \makecell{$-0.001$\\{\scriptsize$(1.59\times10^{-17})$}}\\
        \hline
        \makecell{$\widehat{\sigma}^2_{\mathcal{RE}}$\\{\scriptsize$(SE)$}} & \makecell{$0.684$\\{\scriptsize$(5.14\times10^{-17})$}} & \makecell{$0.741$\\{\scriptsize$(2.71\times10^{-17})$}} & \makecell{$0.850$\\{\scriptsize$(2.58\times10^{-17})$}} & \makecell{$\textbf{0.999}$\\{\scriptsize$(3.94\times10^{-17})$}} & \makecell{$1.214$\\{\scriptsize$(1.37\times10^{-17})$}} & \makecell{$1.359$\\{\scriptsize$(4.74\times10^{-17})$}} & \makecell{$1.423$\\{\scriptsize$(1.71\times10^{-17})$}}\\
        \hline
        \multicolumn{8}{|c|}{For samples with $15\%$ contamination $(0.85g+0.15N(0,16))$.} \\
        \hline
        \makecell{$\widehat{\mu}_{\mathcal{RE}}$\\{\scriptsize$(SE)$}} &\makecell{$-0.005$\\{\scriptsize$(4.38\times10^{-18})$}} & \makecell{$-0.004$\\{\scriptsize$(4.72\times10^{-18})$}}& \makecell{$-0.004$\\{\scriptsize$(1.66\times10^{-17})$}} & \makecell{$-0.004$\\{\scriptsize$(2.22\times10^{-19})$}} & \makecell{$-0.003$\\{\scriptsize$(6.66\times10^{-19})$}} & \makecell{$-0.003$\\{\scriptsize$(1.22\times10^{-18})$}} & \makecell{$-0.002$\\{\scriptsize$(1.77\times10^{-18})$}}\\
        \hline
        \makecell{$\widehat{\sigma}^2_{\mathcal{RE}}$\\{\scriptsize$(SE)$}} & \makecell{$0.823$\\{\scriptsize$(1.59\times10^{-17})$}} & \makecell{$0.892$\\{\scriptsize$(1.66\times10^{-18})$}} &  \makecell{$\textbf{1.022}$\\{\scriptsize$(1.18\times10^{-16})$}} & \makecell{$1.084$\\{\scriptsize$(6.00\times10^{-17})$}} & \makecell{$1.461$\\{\scriptsize$(1.62\times10^{-17})$}} & \makecell{$1.635$\\{\scriptsize$(2.53\times10^{-17})$}} & \makecell{$1.712$\\{\scriptsize$(2.97\times10^{-17})$}} \\
        \hline
        \multicolumn{8}{|c|}{For samples with $20\%$ contamination $(0.80g+0.20N(0,16))$.} \\
        \hline
        \makecell{$\widehat{\mu}_{\mathcal{RE}}$\\{\scriptsize$(SE)$}} & \makecell{$-0.016$\\{\scriptsize$(0)$}}  & \makecell{$-0.016$\\{\scriptsize$(1.02\times10^{-17})$}}&  \makecell{$-0.016$\\{\scriptsize$(7.77\times10^{-18})$}} &  \makecell{$-0.016$\\{\scriptsize$(1.02\times10^{1.19})$}} & \makecell{$-0.015$\\{\scriptsize$(3.66\times10^{-18})$}} & \makecell{$-0.015$\\{\scriptsize$(7.77\times10^{-19})$}} & \makecell{$-0.015$\\{\scriptsize$(3.88\times10^{-19})$}}\\
        \hline
        \makecell{$\widehat{\sigma}^2_{\mathcal{RE}}$\\{\scriptsize$(SE)$}} &  \makecell{$0.929$\\{\scriptsize$(1.18\times10^{-17})$}} & \makecell{$\textbf{1.007}$\\{\scriptsize$(4.77\times10^{-18})$}} & \makecell{$1.155$\\{\scriptsize$(1.41\times10^{-16})$}} & \makecell{$1.359$\\{\scriptsize$(8.28\times10^{-17})$}} & \makecell{$1.654$\\{\scriptsize$(1.87\times10^{-17})$}} & \makecell{$1.852$\\{\scriptsize$(5.30\times10^{-17})$}} & \makecell{$1.941$\\{\scriptsize$(2.05\times10^{-17})$}} \\
        \hline
        \multicolumn{8}{|c|}{For samples with $25\%$ contamination $(0.75g+0.25N(0,16))$.} \\
        \hline
        \makecell{$\widehat{\mu}_{\mathcal{RE}}$\\{\scriptsize$(SE)$}}  & \makecell{$-0.006$\\{\scriptsize$(2.44\times10^{-18})$}} &  \makecell{$-0.006$\\{\scriptsize$(4.33\times10^{-18})$}} &  \makecell{$-0.006$\\{\scriptsize$(1.22\times10^{-18})$}}&   \makecell{$-0.006$\\{\scriptsize$(1.44\times10^{-18})$}} &  \makecell{$-0.004$\\{\scriptsize$(0)$}} &  \makecell{$-0.004$\\{\scriptsize$(3.77\times10^{-18})$}} &  \makecell{$-0.006$\\{\scriptsize$(9.99\times10^{-19})$}}\\
        \hline
        \makecell{$\widehat{\sigma}^2_{\mathcal{RE}}$\\{\scriptsize$(SE)$}}  & \makecell{$\textbf{1.020}$\\{\scriptsize$(4.24\times10^{-17})$}} & \makecell{$1.106$\\{\scriptsize$(1.06\times10^{-16})$}} &  \makecell{$1.270$\\{\scriptsize$(5.98\times10^{-17})$}} & \makecell{$1.496$\\{\scriptsize$(3.85\times10^{-17})$}} & \makecell{$1.825$\\{\scriptsize$(4.68\times10^{-17})$}} &\makecell{$2.048$\\{\scriptsize$(6.15\times10^{-17})$}} & \makecell{$2.147$\\{\scriptsize$(7.77\times10^{-18})$}} \\
        \hline
     \end{longtable}
\begin{table}[ht]
\centering
 \caption{Comparison of different estimators of $(\mu, \sigma^2)$ with different degrees of contamination.}\label{yu}





 \begin{tabular}{|>{\centering\arraybackslash}p{5cm}|>{\centering\arraybackslash}p{5cm}|>{\centering\arraybackslash}p{5cm}|}

 \hline
\textbf{Contamination} & $\mathbf{MDPDE}/\mathbf{MLDPDE}$\hspace{5 in}$\mathbf{(\widehat{\mu}_\mathcal{B}, \widehat{\sigma}_\mathcal{B}^2)}$ &$\mathbf{MLNREE}$ \hspace{5 in}$\mathbf{(\widehat{\mu}_\mathcal{RE}, \widehat{\sigma}_\mathcal{RE} ^2)}$\\
\hline

No contamination & $(-0.012, 0.968)$&$(-0.012, 1.014) $ $\beta = 1.02$ \\

$5\%$ contamination & $(-0.003, 1.141)$&$(-0.003, 1.019)$ $\beta = 0.96$ \\

$10\%$ contamination & $(-0.001, 1.359)$&$(-0.001, 0.999)$ $\beta = 0.91$ \\

$15\%$ contamination & $(-0.003, 1.635)$&$(-0.004, 1.022) $ $\beta = 0.88$ \\

$20\%$ contamination & $(-0.015, 1.852)$&$(-0.016, 1.007)$ $\beta = 0.86$ \\

$25\%$ contamination & $(-0.004, 2.048)$&$(-0.006, 1.020)$ $\beta = 0.85$ \\
 \hline
\end{tabular} \\
\end{table}

\subsection{Student's $r$-distribution with negative degrees of freedom }\label{S7.2}  
Let $Y_1, \dots, Y_n$ be an i.i.d sample from a true distribution $g$, which is a Student's $r$-distribution with the pdf given in Definition \ref{s_d},
where $F_\nu(\sigma^2)=\frac{\Gamma\left(1-\frac{\nu}{2}\right)}{\Gamma\left(1-\frac{\nu+1}{2}\sqrt{-\pi \nu \sigma^2}\right)}$ and $\nu \in (-\infty, 0)$ is assumed to be known. Note that the density can equivalently be written as
\begin{eqnarray*}
  f_\lambda(y) =  \begin{cases}
  N_\nu(\lambda)\left[1+ y^2 \frac{1}{\nu \sigma^2+\mu^2} - y \frac{2 \mu \sigma}{\nu \sigma^2+ \mu^2} \right]^{-\frac{\nu+1}{2}}, & \mu-\sigma\sqrt{-\nu}<y<\mu+\sigma\sqrt{-\nu}\\
     0, & \text{otherwise},
  \end{cases}
 \end{eqnarray*}
where $\lambda=(\mu, \sigma^2)$ and $ N_\nu(\lambda)= \frac{\Gamma\left(1-\frac{\nu}{2}\right)}{\Gamma\left(1-\frac{\nu+1}{2}\sqrt{-\pi \nu \sigma^2}\right)} \left[1+\frac{ \mu^2}{\nu \sigma^2}\right]^{-\frac{2}{\nu+1}}$. Let $\alpha = \beta -\frac{2}{\nu+1}$ and $\frac{1}{\nu} = l_{\alpha, \beta} = \frac{\beta-\alpha}{2+\alpha-\beta}$. Then, the above pdf can equivalently be rewritten as  
\begin{equation}\label{E62}
    f_{\lambda}(y) = \begin{cases}
  Z_{\alpha, \beta}(\lambda) \left[1+ y^2\frac{l_{\alpha, \beta}}{\sigma^2+l_{\alpha,\beta}\mu^2}-y\frac{2\mu\sigma l_{\alpha,\beta}}{\sigma^2+l_{\alpha,\beta}\mu^2}\right]^{\frac{1}{\alpha - \beta}}, & \text{for $\mu -\sigma d_{\alpha, \beta} \leq y \leq \mu + \sigma d_{\alpha, \beta}$}\\
    0, & \text{otherwise}, 
    \end{cases}
\end{equation}
where $Z_{\alpha, \beta} (\lambda)=N_{1/l_{\alpha, \beta}}(\lambda)$ and $d_{\alpha, \beta} = \sqrt{-1/l_{\alpha, \beta}}$. Note that $f_{\lambda}(y) \in \mathcal{M}^{(\alpha, \beta)}$-family (but not regular). Now, we proceed to estimate $\lambda$ by maximizing the generalized likelihood function associated with the LNRE as in \eqref{E16}. 
By using \eqref{E62}, the likelihood function is given by
\begin{eqnarray}
    \mathcal{L}_{\mathcal{RE}}^{(\alpha,\beta)}(y_1 ^n;\lambda)
    &=& \frac{1}{\alpha-\beta}\log \left[\frac{F^{\alpha-\beta}_{1/l_{\alpha, \beta}}(\sigma^2)}{n}\sum_{j=1} ^ {n} \widehat{g}^{\beta-1}(y_j)  \left[1+ l_{\alpha, \beta}\left(\frac{y_j-\mu}{\sigma}\right)^2\right] \textbf{1}(\mu - \sigma d_{\alpha, \beta} \leq y_j \leq \mu + \sigma d_{\alpha, \beta}) \right] \nonumber \\
    &&- \frac{1}{\alpha}\log \left[\sigma F^{\alpha}_{1/l_{\alpha, \beta}}(\sigma^2) \right] - \frac{1}{\alpha}\log \int_{- d_{\alpha, \beta}}^{d_{\alpha, \beta}} [1+ l_{\alpha, \beta} t^2]^{\frac{\alpha}{\alpha-\beta}} dt ,\label{E73}
\end{eqnarray}
where $\textbf{1}(\cdot)$ denotes the indicator function. 
From this, we observe that the second term is independent of $\mu$. Further, maximization of the above likelihood function is equivalent to maximize
\begin{equation}\label{E64}
\ell^{(\alpha,\beta)}(\lambda)  := \frac{{\left[\frac{F^{\alpha-\beta}_{1/l_{\alpha, \beta}}(\sigma^2)}{n}\sum_{j=1} ^ {n} \widehat{g}^{\beta-1}(y_j)  \left[1+ l_{\alpha, \beta}\left(\frac{y_j-\mu}{\sigma}\right)^2\right] \textbf{1}(\mu - \sigma d_{\alpha, \beta} \leq y_j \leq \mu + \sigma d_{\alpha, \beta}) \right]}^{\frac{1}{\alpha-\beta}}}{{\left[\bigintsss_{- d_{\alpha, \beta}}^{d_{\alpha, \beta}} \sigma F^{\alpha}_{1/l_{\alpha, \beta}}(\sigma^2) [1+ l_{\alpha, \beta} t^2]^{\frac{\alpha}{\alpha-\beta}} dt\right]}^{\frac{1}{\alpha}}}.    \end{equation}
Note that $\ell^{(\alpha,\beta)}(\lambda)\neq 0$ if and only if at least one of the $y_j, j =1,\dots, n$, satisfies the condition $\mu - \sigma d_{\alpha, \beta} \leq y_j \leq \mu + \sigma d_{\alpha, \beta}$. 
Thus, $(\widehat{\mu}_{\mathcal{RE}}, \widehat{\sigma}_{\mathcal{RE}}^2)$ can be in any of the regions $\mu - \sigma d_{\alpha, \beta} \leq y_j \leq \mu + \sigma d_{\alpha, \beta}, $ $j=1, \dots ,n$. Now, we calculate the local maximizer in each of these regions in order to find out the global maximizer that corresponds to the MLNREE of $(\mu, \sigma^2)$. Note that the given regions may have nonempty intersections. Thus, we first need to construct disjoint regions and define the likelihood function in each of these regions to find local maximizers and then finally, we can get a global maximizer. However, it is computationally cumbersome task to find disjoint regions and hence, we below consider the problem of estimating one parameter while the other is known. 
\subsubsection{Estimation of $\mu$ when $\sigma^2$ is known}\label{S7.2.1}
We consider the problem of estimating the mean of a Student's $r$-distribution with known $\sigma^2$, say $\sigma^2 =1$. We first draw a sample $Y_1=y_1,\dots,Y_{20}=y_{20}$ from the mixture $0.8g+0.2\mathcal{N}(10,1)$ ($20\%$ contamination), where $g$ is the Student's $r$-distribution with $\nu = -3 \; (\text{i.e., }\alpha - \beta = 1)$, $\mu = 0$ and $\sigma^2 = 1$. Without any loss of generality, let us assume that $y_1\leq \dots\leq y_{20}$.
Now, we have $l_{\alpha, \beta} = \frac{\beta-\alpha}{2+\alpha-\beta}= -\frac{1}{3}$ ($ d_{\alpha, \beta} = \sqrt{3}$) and consequently, the normalizing constant, $F_{1/l_{\alpha, \beta}}$, and the denominator term in the likelihood function given in \eqref{E64} are independent of $\mu$. Thus, for the given estimation problem, it suffices to maximize the likelihood function given by
\begin{equation}\label{as23}
\ell^{(\alpha,\beta)}_1(\mu) := {\sum_{j=1} ^ {n} \widehat{g}^{\beta-1}(y_j)  \left[1- \frac{1}{3} \left({y_j-\mu}\right)^2\right] \textbf{1}(\mu \in I_j) },   
\end{equation}
where $I_j = [y_j - \sqrt{3}, y_j + \sqrt{3}]$, $j \in \{1,\dots, n\}$, and $\widehat{g}$ is the kernel density estimate with Gaussian kernel. Note that the dependency of support upon the parameter $\mu$ results in the appearance of the indicator function in the above likelihood function. 
To maximize this likelihood, we first determine the disjoint intervals and obtain the local maximizers. For $d =1, \dots , n,$ let $I_d' = I_d $\textbackslash$ (\cup_{j=1} ^{d-1} I_j)$. Further, for each $d$, $I_d'$ can be divided into $(n-d+1)$ sub-intervals as $I_d^i:=[(I_d'\cap(\cap_{j=d}^iI_j))$\textbackslash$\cup_{j=i+1}^nI_j],\; i=d,d+1,\dots,n$. For $\mu \in I_d^i$, $i=d, d+1, \dots ,n$, the likelihood function is given by 
$$\ell^{(\alpha,\beta)}_1(\mu) =  \sum_{j=d} ^ {i} \widehat{g}^{\beta-1}(y_j)  \left[1-\frac{1}{3}(y_j-\mu)^2\right], $$
which  is monotonically increasing if $\mu \leq \frac{\sum_{j=d} ^ {i} y_j \widehat{g}^{\beta-1}(y_j)}{\sum_{j=1} ^ {2} \widehat{g}^{\beta-1}(y_j)}$; otherwise, it is monotonically decreasing. Thus, the local maximizer of $\ell^{(\alpha,\beta)}_1(\mu)$ in any nonempty $I_d^i$, $i \in \{d, d+1, \dots ,n\}$, is given by the median of $\bigg\{L_{d}^i, U_d^i, \frac{\sum_{j=d} ^ {i} y_j \widehat{g}^{\beta-1}(y_j)}{\sum_{j=d} ^ {i} \widehat{g}^{\beta-1}(y_j)}\bigg\}$, where $L_d^i$ and $U_d^i$ are the lower and upper bounds of the interval $I_d^i$, respectively. Since $\alpha= 1+\beta > 0$, we are free to choose $\beta \in (-1, \infty)$. Thus, for different choices of the tuning parameter $\beta$, the local maximizer in each interval is obtained and consequently, the global maximizer is derived. For our study, we consider $\beta \in \mathcal{B}$, where $\mathcal{B}=\{0.9,1,1.2,1.3,1.4\}$. A detailed description for obtaining the global maximizer for a single sample of size $20$ from $0.8g+0.2\mathcal{N}(10,1)$ is given in the Appendix-A for illustration. Since we cannot draw a general conclusion by considering only one sample, we replicate this process $1000$ times and report the average values of the MLNREEs along with their standard errors (SE) in Table \ref{T_s1}. Note that the  MLNREEs, for different choices of $\beta \in \mathcal{B}$, and the MDPDE (at $\beta=1$) are very close to the true parameter value  $\mu=0$; the best MLNREE is slightly better than the MDPDE. 
\begin{table}[ht]
\centering
\caption{MLNREEs of $\mu$, for $\beta \in \mathcal{B}$.}
\label{T_s1}
\begin{tabular}{|>{\centering\arraybackslash}p{2.25cm}|>{\centering\arraybackslash}p{2.25cm}|>{\centering\arraybackslash}p{2.25cm}|>{\centering\arraybackslash}p{2.25cm}|>{\centering\arraybackslash}p{2.25cm}|>{\centering\arraybackslash}p{2.25cm}|}
 \hline
$\beta$&$1.4$ & $1.3$ & $1.2$ & $1$ & $0.9$ \\
\hline
\makecell{$\widehat{\mu}_{\mathcal{RE}}$\\{\scriptsize$(SE)$}}&\makecell{$0.016$\\{\scriptsize$(1.24\times10^{-18})$}} & \makecell{$0.016$\\{\scriptsize$(1.27\times10^{-18})$}}& \makecell{$0.016$\\{\scriptsize$(6.66\times10^{-18})$}}& \makecell{$0.017$\\{\scriptsize$(1.00\times10^{-17})$}}& \makecell{$0.017$\\{\scriptsize$(3.19\times10^{-18})$}} \\
\hline
\end{tabular}
\end{table}
\subsubsection{Estimation of $\sigma^2$ when $\mu$ is known:}\label{S7.2.2}
We consider the problem of estimating variance of a Student's $r$-distribution with known mean, say $\mu =0$. We first draw a sample $Y_1=y_1,\dots,Y_{20}=y_{20}$ from the mixture $0.8g+0.2\mathcal{N}(0,25)$, where $g$ is the Student's $r$-distribution with $\nu = -3$  (i.e., $\alpha - \beta = 1$), $\mu = 0$ and $\sigma^2 = 1$. Without any loss of generality, let us assume that $|y_1|\leq \dots \leq|y_{20}|$ (for mathematical simplicity in constructing disjoint intervals). Now, we have $l_{\alpha, \beta} = \frac{\beta-\alpha}{2+\alpha-\beta}= -\frac{1}{3}$ ($ d_{\alpha, \beta} = \sqrt{3}$). Consequently, the likelihood function given in \eqref{E73} reduces to
\begin{eqnarray*}
  \mathcal{L}_{\mathcal{RE}}^{(\alpha,\beta)}(y_1 ^n;\sigma^2) 
    &=&\log \left[\frac{\sqrt{3}}{4 n}\sum_{j=1} ^ {n} \widehat{g}^{\beta-1}(y_j)  \frac{1}{\sqrt{\sigma^2}}\left[1+ l_{\alpha, \beta}\frac{y_j^2}{\sigma^2}\right] \textbf{1}\left(\sigma^2 \geq \frac{y_j^2}{3} \right) \right] \\
    &&- \frac{1}{\alpha}\log (\sigma^2)^{\frac{1-\alpha}{2}} - \frac{1}{\alpha}\log \int_{- \sqrt{3}}^{\sqrt{3}} \left(\frac{\sqrt{3}}{4}\right)^\alpha \left[1- \frac{1}{3}  t^2 \right]^{\alpha} dt,
\end{eqnarray*}
where the third term is independent of $\sigma^2$. Thus, for the given estimation problem, it suffices to consider the likelihood function given by
\begin{eqnarray}
 \ell^{(\alpha,\beta)}_2(\sigma^2) &:=&
  {\sum_{j=1} ^ {n} \widehat{g}^{\beta-1}(y_j) (\sigma^2)^{-\frac{1}{2(1+\beta)}} \left[1- \frac{1}{3} \left(\frac{y_j}{\sigma}\right)^2\right] \textbf{1}(\sigma^2 \in K_j) }\label{E_l},
\end{eqnarray}
where the intervals $K_j = \left[\frac{y_j^2 }{3}, \infty\right)$, $j \in \{1,\dots, n\}$. Next, we can obtain the disjoint intervals, $J_m = \left[\frac{y_m^2}{3}, \frac{y_{m+1}^2}{3}\right]$, $m \in \{1, 2, \dots, 19\}$ and $J_{20}=\left[\frac{y_{20}^2}{3}, \infty\right]$. Now, for $m \in \{1, 2, \dots , 20\}$, let $\sigma^2 \in J_m$. Then, the likelihood function given in \eqref{E_l} reduces to
\begin{equation}\label{E69}
\ell^{(\alpha,\beta)}_2(\sigma^2) = {\sum_{j=1} ^ {m} \widehat{g}^{\beta-1}(y_j) (\sigma^2)^{-\frac{1}{2(1+\beta)}} \left[1-  \frac{y_j^2}{3\sigma^2}\right]}, 
\end{equation}
which is monotonically increasing if 
\begin{equation}\label{E70}
\sigma^2 \leq \frac{3+2\beta}{3}\frac{\sum_{j=1} ^ {m} y_j^2 \widehat{g}^{\beta-1}(y_j)}{\sum_{j=1} ^ {m} \widehat{g}^{\beta-1}(y_j)};
\end{equation}
otherwise, it is monotonically decreasing. Thus, the local maximizer is the median of $$\left\{L_m, U_m, \frac{3+2\beta}{3}\frac{\sum_{j=1} ^ {m} y_j^2 \widehat{g}^{\beta-1}(y_j)}{\sum_{j=1} ^ {m} \widehat{g}^{\beta-1}(y_j)}\right\},$$
where $L_m$ and $U_m$ are the lower and upper bounds of the interval $J_m$, respectively, for $m \in \{1, 2, \dots , 20\}$.
Since $\alpha= 1+\beta > 0$, we are free to choose $\beta \in (-1, \infty)$. Thus, for different choices of the tuning parameter $\beta$, the local maximizer in each interval is obtained and consequently, the global maximizer is derived. For illustration, a detailed description for obtaining the global maximizer for a single sample of size $20$ from $0.8g+0.2\mathcal{N}(0,25)$ is given in Appendix-B. We again replicate the process $1000$ times and report the average values of the MLNREEs along with their standard errors (SE) in Table \ref{T_s2}. 
\begin{table}[ht]
\centering
\caption{MLNREEs of $\mu$, for $\beta \in \mathcal{B}$.}
\label{T_s2}
\begin{tabular}{|>{\centering\arraybackslash}p{2.25cm}|>{\centering\arraybackslash}p{2.25cm}|>{\centering\arraybackslash}p{2.25cm}|>{\centering\arraybackslash}p{2.25cm}|>{\centering\arraybackslash}p{2.25cm}|>{\centering\arraybackslash}p{2.25cm}|}
 \hline
$\beta$&$1.4$ & $1.3$ & $1.2$ & $1$ & $0.9$ \\
\hline
\makecell{$\widehat{\sigma}^2_{\mathcal{RE}}$\\{\scriptsize$(SE)$}}& \makecell{$1.059$\\{\scriptsize$(1.31\times10^{-16})$}} & \makecell{$1.037$\\{\scriptsize$(1.63\times10^{-17})$}} & \makecell{$\textbf{1.012}$\\{\scriptsize$(2.18\times10^{-17})$}} & \makecell{$0.942$\\{\scriptsize$(5.07\times10^{-17})$}} & \makecell{$0.910$\\{\scriptsize$(6.72\times10^{-18})$}} \\
\hline
\end{tabular}
\end{table}
Consequently, the best  MLNREE of $\sigma^2$ is given by $\widehat{\sigma}^2_{\mathcal{RE}}=1.012$, which corresponds to $\beta=1.2$. Further, the MDPDE and MLDPDE of $\sigma^2$ are the same and it is given by $\widehat{\sigma}^2_{\mathcal{B}}=0.942$, which corresponds to $\beta=1$. Thus, the best MLNREE is closer to the true parameter value ($\sigma^2=1$) than the MDPDE/MLDPDE. Hence, the MLNREE is better compared to MDPDE and MLDPDE in presence of outliers, at least for the studied situations.
\begin{table}[H]
\centering
\caption{Newcomb's light speed data}
\label{Table_14} 
\renewcommand{\arraystretch}{1.2} 
\begin{tabular}{|>{\centering\arraybackslash}p{1cm}|>{\centering\arraybackslash}p{1cm}|>{\centering\arraybackslash}p{1cm}|>{\centering\arraybackslash}p{1cm}|>{\centering\arraybackslash}p{1cm}|>{\centering\arraybackslash}p{1cm}|>{\centering\arraybackslash}p{1cm}|>{\centering\arraybackslash}p{1cm}|>{\centering\arraybackslash}p{1cm}|>{\centering\arraybackslash}p{1cm}|>{\centering\arraybackslash}p{1cm}|}
 \hline
 
$26$ & $26$ & $26$  & $26$  & $26$  & $27$  & $27$ & $27$  & $27$  & $27$ & $27$ \\

$25$ & $25$  & $25$  & $25$  & $25$  & $28$  & $28$ & $28$  & $28$  & $28$ & $28$\\

$28$  & $24$  & $24$  & $24$  & $24$  & $24$  & $29$  & $29$  & $29$  & $29$ & $29$\\

$23$  & $23$  & $23$  & $30$  & $30$  & $30$  & $22$ &$22$ & $31$  & $31$ & $21$ \\

 $21$ & $32$  & $32$  & $32$  & $32$  & $32$  & $20$  &$33$  & $33$  & $19$ & $34$\\

$36$ & $36$ &$36$ &$36$ &$16$ &$16$ & $37$ &$39$ &$40$ &$-2$ &$-44$ \\
\hline
\end{tabular}
\end{table}
\subsection{A real data analysis}
\begin{figure} 
\centering
 \includegraphics[width=.8\textwidth]{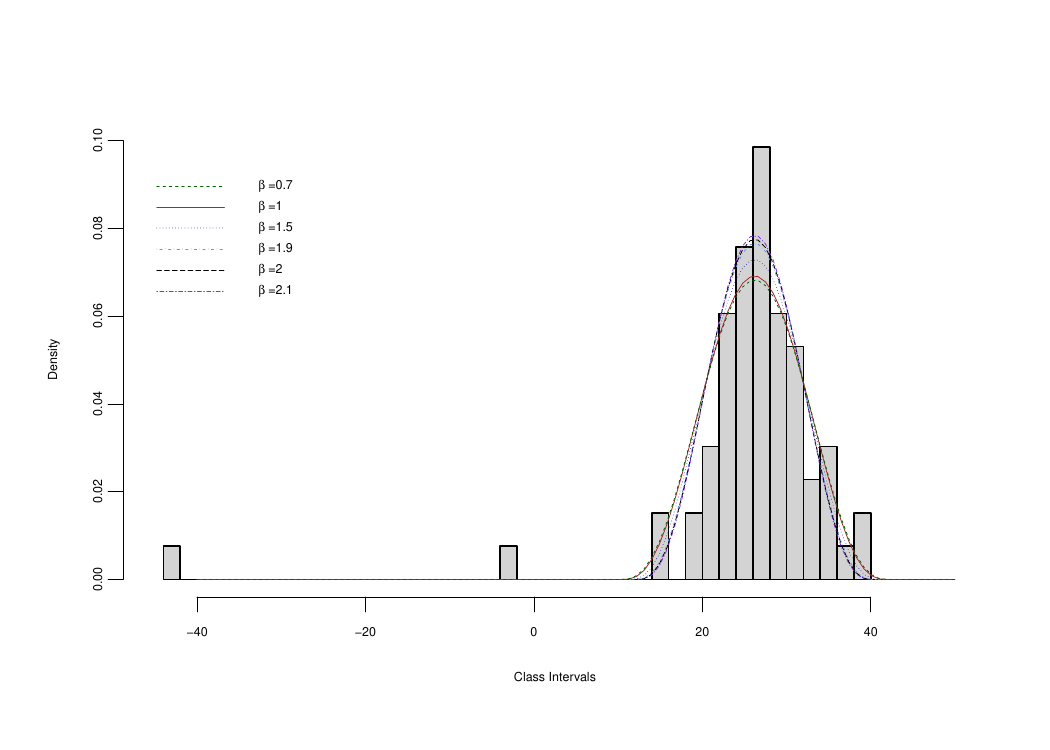}  
    \captionof{figure}{Density plots over the Histogram plot for the Newcomb's data.} 
    \label{Fig_100}
   \end{figure}

To  compare the performance of the MLNREE with that of other related estimators, we consider the Newcomb's light speed data, presented in Table \ref{Table_14}, which is borrowed from \cite{Real_Data}. 
Assume that the underlying true model distribution  $f_\lambda$ is a Student's $r$-distribution with $\nu=-7$. The MLNREE of the parameters $\mu$ and $\sigma^2$ can be obtained by maximizing the likelihood function given in \eqref{E73}. As mentioned in Subsection~\ref{S7.2}, this maximization problem can be solved by constructing disjoint regions, which is a computationally cumbersome task. Thus, 
we consider the problem of estimating $\sigma^2$ when $\mu$ is known. We assume $\mu$ to be the average of the data as it is an unbiased and consistent estimator for the mean of the underlying distribution. Let $y_1, y_2,\dots , y_{66}$ be the given data, as in Table \ref{Table_14}, and consequently, ${\mu} = \overline{y}=26.21$. Now, we have  $\alpha - \beta = 1/3 $ and $l_{\alpha, \beta} = \frac{\beta-\alpha}{2+\alpha-\beta}= -\frac{1}{7}$ ($ d_{\alpha, \beta} = \sqrt{7}$). Further, note that the third term in the likelihood function given in \eqref{E73} is independent of $\sigma^2$. Thus, for the given estimation problem, it suffices to consider the likelihood function given by
\begin{equation*}
 \mathcal{L}_{\mathcal{RE}}^{*(\alpha,\beta)}(y_1 ^n;\sigma^2) :=\log \left((\sigma^2)^{-\frac{3}{2(1+3\beta)}} \left[{\sum_{j=1} ^ {n} \widehat{g}^{\beta-1}(y_j)  \left\{1- \frac{1}{7} \left(\frac{y_j-\widehat{\mu}}{\sigma}\right)^2\right\} \textbf{1}(\sigma^2 \in V_j) } \right]^3\right),    
\end{equation*}
where $V_j = \left[\frac{(y_j-\widehat{\mu})^2 }{7}, \infty\right)$, $j \in \{1,\dots, 66\}$. Maximization of the above likelihood function is equivalent to maximize
\begin{equation}\label{opk}
\ell^{*(\alpha,\beta)}(\sigma^2) :=(\sigma^2)^{-\frac{3}{2(1+3\beta)}}   \left[{\sum_{j=1} ^ {n} \widehat{g}^{\beta-1}(y_j) \left\{1- \frac{1}{7} \left(\frac{y_j-\widehat{\mu}}{\sigma}\right)^2\right\} \textbf{1}(\sigma^2 \in V_j) }\right]^3.   
\end{equation}
 Here, some of the data points are repetitive and consequently, some of the $V_j\text{'}s$ become equal. Further, by proceeding in the same lines as in Subsection~\ref{S7.2}, we construct the set of disjoint intervals $W_m$, $m =1, 2, \dots , 23$, which are given in Table \ref{Table_15}. 
\begin{table}[htp]
\centering
\caption{ Disjoint intervals for finding the MLNREE using the Newcomb's data.}
\label{Table_15}
\renewcommand{\arraystretch}{1.2}  
\begin{tabular}{|>{\centering\arraybackslash}p{4.2cm}|>{\centering\arraybackslash}p{4.2cm}|>{\centering\arraybackslash}p{4.2cm}|>{\centering\arraybackslash}p{4.2cm}|}

\hline
 $W_1=[0.0064, 0.0886]$ & $W_2=[0.0886, 0.2098]$ &$W_3=[0.2098, 0.4566]$ & $W_4=[0.4566, 0.6990]$\\

$W_5=[0.6990, 1.1103]$ & $W_6=[1.1103, 1.4739]$ &$W_7=[1.4739, 2.0497]$ & $W_8=[2.0497, 2.5345]$\\

 $W_9=[2.5345, 3.2748]$ &$W_{10}=[3.2748, 3.8808]$ & $W_{11}[3.8808, 4.7856]$ & $W_{12}=[4.7856, 5.5129]$ \\

$W_{13}=[5.5129, 6.5821]$ & $W_{14}=[6.5821, 7.4306]$ & $W_{15}=[7.4306, 8.6644]$ &$W_{16}=[8.6644, 13.6860]$ \\

$W_{17}=[ 13.6860, 14.8982]$ & $W_{18}=[ 14.8982, 16.6254]$ & $W_{19}=[16.6254, 23.3614]$ & $W_{20}=[23.3614, 27.1579] $\\

 $W_{21}=[27.1579, 113.703]$ & $W_{22}=[113.703, 704.248]$ & $W_{23}=[704.2489, \infty)$ & \\
\hline
\end{tabular}
\end{table}
Let $\sigma^2 \in W_1$. Then,  the likelihood function given in \eqref{opk} can be written as
\begin{equation}\label{E71}
\ell^{*(\alpha,\beta)}(\sigma^2) = (\sigma^2)^{-\frac{3}{2(1+3\beta)}}\left[{\sum_{j=1} ^ {5} \widehat{g}^{\beta-1}(y_j) \left\{1- \frac{1}{7} \left(\frac{y_j-\widehat{\mu}}{\sigma}\right)^2\right\}}\right]^3,
\end{equation}
which is monotonically increasing if 
\begin{equation}\label{E72}
\sigma^2 \leq \frac{3+6\beta}{7}\frac{\sum_{j=1} ^ {5} (y_j-\widehat{\mu})^2 \widehat{g}^{\beta-1}(y_j)}{\sum_{j=1} ^ {5} \widehat{g}^{\beta-1}(y_j)};
\end{equation}
otherwise, it is monotonically decreasing. Thus, the local maximizer in the interval $W_1$ is the median of $$\left\{L_1, U_1,\frac{3+6\beta}{7}\frac{\sum_{j=1} ^ {5} (y_j-\widehat{\mu})^2 \widehat{g}^{\beta-1}(y_j)}{\sum_{j=1} ^ {5} \widehat{g}^{\beta-1}(y_j)}\right\},$$
where $L_1$ and $U_1$ are the lower and upper bounds of $W_1$, respectively. By proceeding in the same lines, we get the local maximizer for each disjoint interval given in Table~\ref{Table_15}. Note that these local maximizers depend on $\beta$. Since $\alpha= \frac{1}{3}+\beta > 0$, we are free to choose $\beta \in (-1/3, \infty)$ and consequently, $\alpha$ is determined from the relation $\alpha= \frac{1}{3}+\beta$. For our study, we consider $\beta \in \mathcal{C}$, where $\mathcal{B}=\{0.9, 1, 1.5, 1.9, 2, 2.1\}$. Now, from
the set of all local maximizers obtained corresponding to $23$ disjoint sub-intervals, we get the the global
maximizer, as in Table \ref{Table_18} , for each $\beta \in \mathcal{C}$, which is indeed the MLNREE of $\sigma^2$, for $\beta \in \mathcal{C}$. Further, the MDPDEs and MLDPDEs of $\sigma^2$ are the same and it is given by $\widehat{\sigma}_{\mathcal{B}}^2=35.74$, which corresponds to $\beta=1$.
\begin{table}[ht]
\centering
\caption{MLNREEs of $\sigma^2$, for $\beta \in \mathcal{C}$ based on the Newcomb's data.}
\label{Table_18} 
\begin{tabular}{|>{\centering\arraybackslash}p{2cm}|>{\centering\arraybackslash}p{2.2cm}|>{\centering\arraybackslash}p{2.2cm}|>{\centering\arraybackslash}p{2.2cm}|>{\centering\arraybackslash}p{2.2cm}|>{\centering\arraybackslash}p{2.2cm}|>{\centering\arraybackslash}p{2.2cm}|}
 \hline

$\beta$&$0.9$ & $1$ & $1.5$ & $1.9$ & $2$ & $2.1$\\
\hline
$\widehat{\sigma}^2_{\mathcal{RE}}$&$36.23$ & $35.74$ & $32.13$ & $29.12$ & $28.43$ & $27.78$\\
\hline
$D_{KS}$& $0.1540$ & $0.1530$ & $0.1451$& $0.1375$ & $\mathbf{0.1366}$ & $0.1373$\\
\hline
\end{tabular}
\end{table}
To see the performance of the MLNREE, we plot the pdf of the given Student's $r$-distribution for different estimated parameter values that are obtained corresponding to different $\beta$ values. From Figure \ref{Fig_100}, we observe that, for the MDPDE and MLDPDE of $\sigma^2$, the density curve shows a higher variability compared to the variability that is there in the histogram plot of the data. On the other hand, the MLNREEs are obtained for different values of $\beta$. Consequently, we have the flexibility to choose a $\beta$ value so that the distribution is fitted well to the given data. To draw a conclusion about the best fitted model to the given data,  we calculate the Kolmogorov-Smirnov (KS) distance given by $D_{KS}=\sup_x|F_n(x)-F(x)|$,  for each MLNREE given in Table \ref{Table_18}, where $F_n$ is the empirical cumulative distribution function (CDF) of the sample and $F$ is the CDF of the fitted distribution. Figure~\ref{Fig_100} shows that, for $\beta \in \{1.9,2,2.1\}$, the Student's $r$-distribution with the estimated $\widehat{\sigma}^2_{\mathcal{RE}}$ fits the given  data quite well. However, the MLNREE for $\beta=2$ outperforms other MLNREEs as it has the minimum KS-distance. Thus, we may conclude that the best MLNREE of $\sigma^2$ is given by $\widehat{\sigma}^2_{\mathcal{RE}}=28.43$.
\section{Concluding remarks}\label{S8}
Minimum divergence methods are widely used in robust statistical inference. 
  The estimation of parameters by minimizing a divergence is equivalent to maximization of the associated generalized likelihood function. For example, the popularly used maximum likelihood estimator is obtained by minimizing the KLD or equivalently, by maximizing the usual log-likelihood function. In this paper, we consider the estimation problems based on a generalized likelihood function associated with the LNRE. Having an additional tuning parameter, the minimum LNRE estimation usually gives more robust parameter estimate compared to those of the minimum DPD and minimum LDPD estimation for data contaminated  by outliers. In particular, we study some sufficiency principles associated with the minimum LNRE estimation. We first introduce  a two-parameter power-law family of distributions, which we call the  $\mathcal{M}^{(\alpha,\beta)}$-family. This family contains the usual exponential family, the $\mathcal{B}^{(\alpha)}$-family, and the $\mathcal{M}^{(\alpha)}$-family as special cases. Then, we prove the Fisher-Darmois-Koopman-Pitman theorem which states that the $\mathcal{M}^{(\alpha,\beta)}$-family has a fixed number of sufficient statistics, independent of sample size, with respect to the generalized likelihood function associated with the LNRE. Further, we obtain the generalized  minimal sufficient statistic for regular $\mathcal{M}^{(\alpha,\beta)}$-family. Then we derive the generalized Rao-Blackwell theorem which is useful to obtain the best estimator for $\mathcal{M}^{(\alpha,\beta)}$-family. Further, we establish the generalized Cram\'{e}r-Rao lower bound for the minimum LNRE estimation. We show that the generalized minimal sufficient statistic for regular $\mathcal{M}^{(\alpha,\beta)}$-family does not attain this bound, apart from a particular case. 
  Further, we obtain the MLNREEs for the family of Student's $t$-distributions . Then, we conduct simulation studies to show that the MLNREEs for the family of Student's $t$ as well as Student's $r$-distributions are more robust compared to the MDPDEs and MLDPDEs. At the end, we present a real data analysis, further illustrating the effectiveness of the minimum LNRE estimation in model building for data contaminated by outliers.
\\\hspace*{0.2 in}Apart from the study done in this paper, there are many potential problems, for the minimum LNRE estimation, which are yet to be explored. For example, asymptotic properties of the MLNREE were studied in the literature for discrete case only. The same for continuous case has not yet been derived in the literature. We are currently working on this problem and hope to report our findings in a future communication.
 \\\\{\bf Acknowledgments}\\
\hspace*{0.2 in}The first author sincerely acknowledges the financial support received from the Ministry of Education, Government of India, under the Prime Minister's Research Fellowship (PMRF) scheme. 
\\\\{\bf Conflicts of Interest}\\
\hspace*{0.2 in} The authors declare no conflict of interest.
 
\noindent \textbf{Appendix-A:}
\\\\Here we give a detailed description of the estimation problem as considered in Subsection \ref{S7.2.1} for a single random sample of size $20$. Consider the sample given in Table \ref{Table_3}. We first construct the intervals $I_j=[y_j-\sqrt{3}, y_j+\sqrt{3}]$, $j=1,2,\dots,20$, 
\begin{table}[ht]
\centering
\caption{ A sample from the Student's $r$-distribution with $\nu = -3 $, $\mu = 0$ and $\sigma = 1$.}
\label{Table_3} 
\renewcommand{\arraystretch}{1.2} 
\begin{tabular}{|>{\centering\arraybackslash}p{3.4cm}|>{\centering\arraybackslash}p{3.4cm}|>{\centering\arraybackslash}p{3.4cm}|>{\centering\arraybackslash}p{3.4cm}|}
 \hline
$-1.7827$ & $-1.1761$ &$-1.0597$ & $-0.3236$\\

$-0.2340$ & $0.4706$ &$0.4712$ & $0.5435$\\

$0.6309$ & $0.7533$ &$0.8020$ & $0.9237$\\

$1.1394$ & $1.4373$ &$1.5351$ & $1.6941$\\

$8.6501$ & $10.7254$ &$12.7694$ & $13.0349$\\
\hline
\end{tabular}
\end{table}
as in Table \ref{Table_4}, which shows that
\begin{table}[ht]
\centering
\caption{ Intervals $I_j = [y_j - \sqrt{3}, y_j + \sqrt{3}]$, $j \in \{1,\dots, 20\}$ for the data given in Table \ref{Table_3}.}
\label{Table_4}
\renewcommand{\arraystretch}{1.2} 
\begin{tabular}{|>{\centering\arraybackslash}p{4cm}|>{\centering\arraybackslash}p{4cm}|>{\centering\arraybackslash}p{4cm}|>{\centering\arraybackslash}p{4cm}|}
 \hline
$I_1=[-3.4607, 0.0033]$ & $I_2=[-2.9081, 0.5559]$ &$I_3=[-2.7917, 0.6723]$ & $I_4=[-2.0556, 1.4084]$\\

$I_5=[-1.9660, 1.4980]$ & $I_6=[-1.2614, 2.2026]$ &$I_7=[-1.2608, 2.2032]$ & $I_8=[-1.1885, 2.2755]$\\

$I_9=[-1.1011, 2.3629]$ & $I_{10}=[-0.9787, 2.4853]$ &$I_{11}=[-0.9300, 2.5340]$ & $I_{12}=[-0.8083, 2.6557]$\\

$I_{13}=[-0.5926, 2.8714]$ & $I_{14}=[-0.2947, 3.1693]$ &$I_{15}=[-0.1969, 3.2671]$ & $I_{16}=[-0.0379, 3.4261]$\\

$I_{17}=[6.9180, 10.3821]$ & $I_{18}=[8.9933, 12.4574]$ &$I_{19}=[11.0373, 14.5014]$ & $I_{20}=[11.3028, 14.7669]$\\
\hline
\end{tabular}
\end{table}
 the first $16$ intervals, $I_1,\dots,I_{16}$, have nonempty intersection as well as  the last four intervals, $I_{17},\dots,I_{20}$, have nonempty intersection. Next, we construct disjoint intervals using the formulation given in Subsection \ref{S7.2.1}. Note that $I_1$ can be divided into $16$ disjoint intervals, namely, $I_1^1,\dots,I_1^{16}$. Similarly, for $I_2$, the sub-interval is given by $I_2'=I_2^{16}=[0.0033,0.5559]$ as other sub-intervals are empty. By proceeding in the same lines, we construct the disjoint sub-intervals for each of $I_3,\dots,I_{20}$. Overall, we get  $38$ disjoint intervals as given in Table \ref{Table_5}. 
\begin{table}[ht]
\centering
\caption{Disjoint intervals for the data given in Table \ref{Table_3}.} 
\label{Table_5} 
\renewcommand{\arraystretch}{1.5} 
\begin{tabular}{|>{\centering\arraybackslash}p{4.2cm}|>{\centering\arraybackslash}p{4.2cm}|>{\centering\arraybackslash}p{4.2cm}|>{\centering\arraybackslash}p{4.2cm}|}
 \hline 

 $I_1^1=[-3.4607, -2.9081]$ & $I_1^2=[-2.9081,-2.7917]$ & $I_1^3=[-2.7917, -2.0556]$ & $I_1^4=[-2.0556, -1.9660]$\\
 
$I_1^5=[-1.9660, -1.2614]$& $I_1^6=[-1.2614, -1.2608]$ & $I_1^7=[-1.2608, -1.1885]$  & $I_1^8=[-1.1885, -1.1011]$\\

$I_1^9=[-1.1011, -0.9787]$  & $I_1^{10}=[-0.9787, -0.9300]$ & $I_1^{11}=[-0.9300, -0.8083]$ & $I_1^{12}=[-0.8083, -0.5926]$\\

$I_1^{13}=[-0.5926, -0.2947]$& $I_1^{14}=[-0.2947, -0.1969]$  & $I_1^{15}=[-0.1969, -0.0379]$ & $I_1^{16}=[-0.0379, 0.0033]$\\

$I_2^{16}=[0.0033, 0.5559]$  & $I_3^{16}=[0.5559, 0.6723]$ &$I_4^{16}=[0.6723, 1.4084]$ &$I_5^{16}=[1.4084, 1.4980]$ \\

$I_6^{16}=[1.4980, 2.2026]$ & $I_7^{16}=[2.2026, 2.2032]$ & $I_8^{16}=[2.2032, 2.2755]$& $I_9^{16}=[2.2755, 2.3629]$\\

 $I_{10}^{16}=[2.3629, 2.4853]$& $I_{11}^{16}=[2.4853, 2.5340]$ & $I_{12}^{16}=[2.5340, 2.6557]$ & $I_{13}^{16}=[2.6557, 2.8714]$\\
 
$I_{14}^{16}=[2.8714, 3.1693]$ & $I_{15}^{16}=[3.1693, 3.2671]$ & $I_{16}^{16}=[3.2671, 3.4261]$ & $I_{17}^{17}=[6.9180, 8.9933]$\\

$I_{17}^{18}=[8.9933, 10.3821]$& $I_{18}^{18}=[10.3821, 11.0373]$ & $I_{18}^{19}=[11.0373, 11.3028]$& $I_{18}^{20}=[11.3028, 12.4574]$\\

$I_{19}^{20}=[12.4574, 14.5014]$ & $I_{20}^{20}=[14.5014, 14.7669]$ & & \\
\hline
\end{tabular}
\end{table}
 Let  $\mu \in I_1^1 = [-3.4607, -2.9081]$. Then, the local maximizer is the median of $\{-3.4607,-2.9081,y_1\}$. If $\mu \in I_1^2$ then the local maximizer is the median of 
$$\Bigg\{ -2.9081, -2.7917, \frac{\sum_{j=1} ^ {2} y_j \widehat{g}^{\beta-1}(y_j)}{\sum_{j=1} ^ {2} \widehat{g}^{\beta-1}(y_j)}\Bigg\}.$$ By proceeding similarly, we get the local maximizer for each interval given in Table \ref{Table_5}. Note that these maximizers depend on $\beta$. Since $\alpha= 1+\beta > 0$, we are free to choose $\beta \in (-1, \infty)$ and consequently, $\alpha$  is determined from the relation $\alpha= 1+\beta$. Now, from the set of all local maximizers obtained corresponding to $38$ disjoint sub intervals,
we get the global maximizer, as in Table~\ref{Table_7}, for different choices of $\beta$. 
\begin{table}[ht]
\centering
\caption{MLNREEs of $\mu$, for different choices of $\beta$ for the data given in Table \ref{Table_3}.}
\label{Table_7} 
\begin{tabular}{|>{\centering\arraybackslash}p{2.25cm}|>{\centering\arraybackslash}p{2.25cm}|>{\centering\arraybackslash}p{2.25cm}|>{\centering\arraybackslash}p{2.25cm}|>{\centering\arraybackslash}p{2.25cm}|>{\centering\arraybackslash}p{2.25cm}|}
\hline
$\beta$&$1.4$ & $1.3$ & $1.2$ &$1$ & $0.9$ \\
\hline
$\widehat{\mu}_\mathcal{RE}$&$0.767$ & $0.765$ & $0.763$ & $0.757$ & $0.753$\\
\hline
\end{tabular}
\end{table}

\vspace*{0.2 in}
\noindent \textbf{Appendix-B:}
\\\\Here we give a detailed description of the estimation problem as considered in Subsection \ref{S7.2.2} for a single random sample of size $20$. Consider the sample given in Table \ref{Table_8}.
\begin{table}[ht]
\centering
\caption{A sample from the Student's $r$-distribution with $\nu = -3$, $\mu = 0$ and $\sigma^2 = 1$.}
\label{Table_8} 
\renewcommand{\arraystretch}{1.2} 
\begin{tabular}{|>{\centering\arraybackslash}p{3.4cm}|>{\centering\arraybackslash}p{3.4cm}|>{\centering\arraybackslash}p{3.4cm}|>{\centering\arraybackslash}p{3.4cm}|}
 \hline
$0.0168$ & $0.0593$ & $-0.1015$ & $-0.4325$ \\

$-0.4620$ &$-0.4669$ & $-0.5620$ & $0.6270$ \\

$-0.6851$ &$-0.7338$ & $0.7675$ & $-0.7915$\\

$0.8283$ &$-0.8574$ & $1.1092$ & $-1.2844$\\

$-1.3149$ &$1.8781$ & $-3.2350$ & $4.3409$\\
\hline
\end{tabular}
\end{table}
 We first construct the intervals $K_j = \left[\frac{y_j^2 }{3}, \infty\right)$, $j \in \{1,\dots, n\}$, as in Table \ref{Table_9}. Next, we obtain the disjoint intervals, $J_m$, $m \in \{1, 2, \dots, 20\}$, as in Table \ref{Table_10}. Let $\sigma^2 \in J_1$. Then, the local maximizer is the median of $\left\{0.00009, 0.00117, \frac{(3+2\beta)y_1^2}{3}\right\}$. If $\sigma^2 \in J_2$, then the local maximizer is the median of $$\left\{0.00117, 0.00343, \frac{3+2\beta}{3}\frac{\sum_{j=1}^{2} y_j^2 \widehat{g}^{\beta-1}(y_j)}{\sum_{j=1}^{2} \widehat{g}^{\beta-1}(y_j)}\right\}.$$
By proceeding similarly, we get local maximizer for each interval  given in Table \ref{Table_10}. Note that these maximizers depend on $\beta$. Since $\alpha= 1+\beta > 0$, we are free to choose $\beta \in (-1, \infty)$ and consequently, $\alpha$  is determined from the relation $\alpha= 1+\beta$. Now, from the set of all local maximizers obtained corresponding to $20$ disjoint sub intervals, we get the global maximizer, as in Table~\ref{Table_13}, for different choices of $\beta$. 
\begin{table}[ht]
\centering
\caption{ Intervals $K_j = \left[\frac{y_j^2 }{3}, \infty\right)$, $j \in \{1,\dots, 20\}$ for the data given in Table \ref{Table_8}. }
\label{Table_9}
\renewcommand{\arraystretch}{1.2} 
\begin{tabular}{|>{\centering\arraybackslash}p{4cm}|>{\centering\arraybackslash}p{4cm}|>{\centering\arraybackslash}p{4cm}|>{\centering\arraybackslash}p{4cm}|}
 \hline
$K_1=[0.00009, \infty)$ & $K_2=[0.00117, \infty)$ &$K_3=[0.00343, \infty)$ & $K_4=[0.06235, \infty)$\\

$K_5=[0.07116, \infty)$ & $K_6=[0.07268, \infty)$ &$K_7=[0.1052, \infty)$ & $K_8=[0.1310, \infty)$\\

$K_9=[0.1564, \infty)$ & $K_{10}=[0.17951, \infty)$ &$K_{11}=[0.1963, \infty)$ & $K_{12}=[0.2088, \infty)$\\

 $K_{13}=[0.2287, \infty)$ &$K_{14}=[0.2450, \infty)$ & $K_{15}=[0.4101, \infty)$ & $K_{16}=[0.5499, \infty)$\\

$K_{17}=[0.5764, \infty)$ &$K_{18}=[1.1758, \infty)$ & $K_{19}=[3.4886, \infty)$ & $K_{20}=[6.2811, \infty)$\\
\hline
\end{tabular}
\end{table}
\begin{table}[ht]
\centering
\caption{ Disjoint intervals for the data given in Table \ref{Table_8}.}
\label{Table_10}
\renewcommand{\arraystretch}{1.2} 
\begin{tabular}{|>{\centering\arraybackslash}p{4cm}|>{\centering\arraybackslash}p{4cm}|>{\centering\arraybackslash}p{4cm}|>{\centering\arraybackslash}p{4cm}|}
 \hline
$J_1=[0.00009, 0.00117]$ & $J_2=[0.00117, 0.00343]$ &$J_3=[0.00343, 0.06235]$ & $J_4=[0.06235, 0.07116]$\\

$J_5=[0.07116, 0.07268]$ & $J_6=[0.07268, 0.1052]$ &$J_7=[0.1052, 0.1310]$ & $J_8=[0.1310, 0.1564]$\\

$J_9=[0.1564, 0.17951]$ & $J_{10}=[0.17951, 0.1963]$ &$J_{11}=[0.1963, 0.2088]$ & $J_{12}=[0.2088, 0.2287]$\\

 $J_{13}=[0.2287, 0.2450]$ &$J_{14}=[0.2450, 0.4101]$ & $J_{15}=[0.4101, 0.5499]$ & $J_{16}=[0.5499, 0.5764]$\\

$J_{17}=[0.5764, 1.1758]$ &$J_{18}=[1.1758, 3.4886]$ & $J_{19}=[3.4886, 6.2811]$ & $J_{20}=[6.2811, \infty)$\\
\hline
\end{tabular}
\end{table}
\begin{table}[H]
\centering
\caption{MLNREEs of $\sigma^2$, for different choices of $\beta$ for the data given in Table \ref{Table_8}.}
\label{Table_13} 
\begin{tabular}{|>{\centering\arraybackslash}p{2.25cm}|>{\centering\arraybackslash}p{2.25cm}|>{\centering\arraybackslash}p{2.25cm}|>{\centering\arraybackslash}p{2.25cm}|>{\centering\arraybackslash}p{2.25cm}|>{\centering\arraybackslash}p{2.25cm}|}
 \hline
$\beta$&$1.4$ & $1.3$ & $1.2$ & $1$ & $0.9$ \\
\hline
$\widehat{\sigma}^2_{\mathcal{RE}}$&$1.0217$ & $\textbf{1.0026}$ & $0.9828$ & $0.9408$ & $0.9184$ \\
\hline
\end{tabular}
\end{table}
\end{document}